\documentclass[review,onefignum,onetabnum]{siamonline220329}
\usepackage{amsmath, amssymb, latexsym, amscd, amsfonts,amstext}
\usepackage[mathscr]{eucal}
\usepackage{graphicx}
\usepackage{subfig}
\usepackage{lipsum}
\usepackage{amsfonts}
\usepackage{graphicx}
\usepackage{epstopdf}
\usepackage{algorithmic}
\usepackage{lipsum}
\usepackage{amsfonts}
\usepackage{graphicx}
\usepackage{epstopdf}
\usepackage{algorithmic}
\usepackage{enumitem}
\usepackage{amsopn}

\ifpdf
\DeclareGraphicsExtensions{.eps,.pdf,.png,.jpg}
\else
\DeclareGraphicsExtensions{.eps}
\fi

\numberwithin{equation}{section}

\allowdisplaybreaks[4]
\ifpdf
\DeclareGraphicsExtensions{.eps,.pdf,.png,.jpg}
\else
\DeclareGraphicsExtensions{.eps}
\fi
\setlist[enumerate]{leftmargin=.5in}
\setlist[itemize]{leftmargin=.5in}

\headers{Convexification method for CIP for RRTE}{M. V. Klibanov, J. Li, L. H. Nguyen,  V. G. Romanov and Z. Yang}

\ifpdf
\hypersetup{
pdftitle={Convexification method for CIP for RRTE},
pdfauthor={M. V. Klibanov, J. Li, L. H. Nguyen,  V. G. Romanov and Z. Yang}
}
\fi
\input{tcilatex}
\begin{document}

\nolinenumbers 

\title{Convexification Numerical Method for a Coefficient Inverse Problem
for the Riemannian Radiative Transfer Equation \thanks{%
2023.04.12 
\funding{The work of J. Li was partially supported by the NSF of China No. 11971221, Guangdong NSF Major Fund No. 2021ZDZX1001, the Shenzhen Sci-Tech Fund No. RCJC20200714114556020, JCYJ20200109115422828 and JCYJ20190809150413261.	
The work of L.H. Nguyen was partially supported by National Science Foundation grant DMS-2208159 and by funds provided by the Faculty Research Grant program at University of North Carolina at Charlotte, Fund No. 111272.
The work by V.G. Romanov was performed within the state assignment of the Sobolev Institute of Mathematics of the Siberian Branch of the Russian Academy of Science, project number FWNF-2022-0009.}%
}}
\author{ Michael V. Klibanov \thanks{
	Department of Mathematics and Statistics, University of North Carolina at
	Charlotte, Charlotte, NC, 28223, USA (\email{mklibanv@uncc.edu}, \email{loc.nguyen@uncc.edu})}
\and Jingzhi Li \thanks{
	Department of Mathematics \& National Center for Applied Mathematics
	Shenzhen \& SUSTech International Center for Mathematics, Southern
	University of Science and Technology, Shenzhen 518055, P.~R.~China (\email{li.jz@sustech.edu.cn})}
\and Loc H. Nguyen \footnotemark[2]
\and Vladimir G. Romanov \thanks{
	Sobolev Institute of Mathematics, Novosibirsk, 630090, Russian Federation (\email{romanov@math.nsc.ru})}
\and Zhipeng Yang \thanks{
	Department of Mathematics, Southern University of Science and Technology,
	Shenzhen 518055, P.~R.~China (\email{yangzp@sustech.edu.cn})} }
\maketitle

\begin{abstract}
The first globally convergent numerical method for a Coefficient Inverse
Problem (CIP) for the Riemannian Radiative Transfer Equation (RRTE) is
constructed. This is a version of the so-called \textquotedblleft
convexification" method, which has been pursued by this research group for a
number of years for some other CIPs for PDEs. Those PDEs are significantly
different from the RRTE. The presence of the Carleman Weight Function (CWF)
in the numerical scheme is the key element which insures the global
convergence. Convergence analysis is presented along with the results of
numerical experiments, which confirm the theory. RRTE governs the
propagation of photons in the diffuse medium in the case when they propagate
along geodesic lines between their collisions. Geodesic lines are generated
by the spatially variable dielectric constant of the medium.
\end{abstract}




\section{Introduction}

\label{sec:1}

The conventional steady state radiative transfer equation (RTE) governs
light propagation in the diffuse medium, such as, e.g. turbulent atmosphere
and biological media \cite{Heino}. Inverse problems for the RTE have
applications in, e.g. problems of seeing through a turbulent atmosphere and
in an early medical diagnostics. In the latter case the near infrared light
with a relatively small energy of photons is used, see, e.g. \cite{Bal}.
However, it is assumed in the RTE that photons propagate along straight
lines between their collisions. On the other hand, since the dielectric
constants in heterogeneous media, such as, e.g. ones mentioned above, vary
in space, then photons actually propagate along geodesic lines between their
collisions. These lines are generated by the Riemannian metric $\sqrt{%
\varepsilon _{r}\left( \mathbf{x}\right) }\left\vert d\mathbf{x}\right\vert
. $ Here and below $\mathbf{x}=\left( x,y,z\right) \in \mathbb{R}^{3}$ and $%
\varepsilon _{r}\left( \mathbf{x}\right) $ is the spatially distributed
dielectric constant, so \ that $n\left( \mathbf{x}\right) =\sqrt{\varepsilon
_{r}\left( \mathbf{x}\right) }$ is the refractive index. To take this into
account, the so-called Riemannian Radiative Transfer Equation (RRTE) should
be used.

This is the first publication, in which a globally convergent numerical
method, the so-called convexification method, is constructed for a
Coefficient Inverse Problem (CIP) for the steady state RRTE. In the past,
numerical methods for inverse problems for the steady state RTE were mostly
developed for the case of inverse source problems \cite{HT1,HT2,Smirnov}.
Inverse source problems are linear. On the other hand, CIPs are nonlinear.
We refer to two recent publications of this research team \cite{KTR,KTR1}
for two versions of the convexification numerical method for a CIP for the
RTE. The presence of the Riemannian aspect in the RRTE causes significant
additional difficulties for the corresponding CIP, as compared with the case
of the RTE in \cite{KTR,KTR1}. The authors are unaware about other numerical
methods for CIPs neither for the RTE nor for the RRTE.

Various uniqueness and stability results for inverse problems for both RTE
and RRTE, including quite general forms of the latter equation, were
published in the past. Since this paper is concerned only with a numerical
method, then we refer now only to a limited number of such publications \cite%
{BalReview,Bal1,Bal,GY,KY,KP,Lay,McD}.

The phenomena of ill-posedness and nonlinearity of CIPs are well known and
cause serious challenges for their numerical solutions. Both a powerful and
popular concept of numerical methods for CIPs\ is based on the minimization
of appropriate least squares cost functionals, see, e.g. \cite%
{B2,B3,B4,Gonch1,Gonch2,Giorgi,Hassi} and references cited therein. Since
such a cost functional is typically non convex, then it usually suffers from
the phenomenon of local minima and ravines, see, e.g. \cite{Scales}, i.e.
the availability of a good first guess about the true solution is a
necessary assumption of the convergence analysis of these numerical methods.

\textbf{Remark 1.1}. \emph{We call a numerical method for a CIP globally
convergent if a theorem is proven, which claims that this method delivers at
least one point in a sufficiently small neighborhood of the true solution
without any advanced knowledge of this neighborhood. The size of that
neighborhood should depend only on the level of noise in the data. }

The key element of our numerical method is the presence of a Carleman Weight
Function (CWF) in a certain weighted least squares cost functional. This
presence ensures the global strict convexity of that functional. This is why
we call our method \textquotedblleft convexification". The CWF is the
function, which is involved as the weight function in the Carleman estimate
for the corresponding PDE operator. Our convergence analysis ensures the
global convergence of the gradient descent method of the minimization of
that functional to the true solution of our CIP, as long as the level of the
noise in the data tends to zero. The apparatus of the Riemannian geometry is
also used here. Results of numerical experiments are presented, and they
confirm our theory.

The convexification concept generates globally convergent numerical methods
since these methods do not rely on good first guesses about the solutions.
The convexification was originally proposed in purely theoretical works \cite%
{KlibIous, Klib97}. Its active numerical studies have started in 2017 after
the publication \cite{Bak}, which has removed some obstacles for numerical
implementations. In this regard, we refer to, e.g. \cite{Khoa,KL,KTR,KTR1}
and references cited therein.

Another important new element of this paper is Theorem 1 (section 3), which
claims existence, uniqueness and positivity of the solution of for the
forward problem for the RRTE. An analog of this theorem for the
non-Riemannian case was proven in \cite{KTR}. The proof of Theorem 1 is
constructive since it ends up with an analysis of a linear integral equation
of the Volterra type. This equation is quite helpful in our numerical
studies in section 6, since we solve it numerically to computationally
simulate the data for the inverse problem. It is well known that such
computational simulations form an important part of numerical studies of any
inverse problem. The presence of the Riemannian aspect creates a significant
additional difficulty in the proof of Theorem 1, as compared with the case
of RTE in \cite{KTR}. This difficulty is due to the necessity of working
with the differential geometry, which, however, was not necessary to do in 
\cite{KTR}.

As to the apparatus of Carleman estimates, it was introduced in the field of
CIPs in the publication \cite{BukhKlib}, initially with the single goal of
proofs of uniqueness theorems. Since then the idea of \cite{BukhKlib} was
explored in many other publications, see, e.g. \cite%
{BY,Crist,Fu,GY,Is,Klib92,Ksurvey,KL,Lay,Yam} and references cited therein.
The convexification principle represents an extension of the idea of \cite%
{BukhKlib} to the topic of globally convergent numerical methods for CIPs.
Those numerical methods might be generalized and employed for important
applications like, e.g. cloaking and quantum scattering studied in \cite%
{LLRU15, LLM21}.

We consider below only real valued functions. For the sake of definiteness,
we work below in our theoretical derivations only with the 3d case. On the
other hand, we present numerical results in the 2d case since the theory for
the 2d case is completely similar with the one in the 3d case. In section 2
we pose both the forward and inverse problems for RRTE. In section 3 we
formulate and prove the above mentioned Theorem 1. In section 4 we derive a
version of the convexification method for our CIP. In section 5 we provide
convergence analysis. Section 6 is devoted to numerical studies, which
confirm our theory.

\section{Statements of Forward and Inverse Problems}

\label{sec:2}

Let numbers $A,a,b,d>0$, where 
\begin{equation}
0<a<b.  \label{2.1}
\end{equation}%
Define the rectangular prism $\Omega \subset \mathbb{R}^{3}$ and parts $%
\partial _{1}\Omega ,\partial _{2}\Omega ,\partial _{3}\Omega $ of its
boundary $\partial \Omega ,$ as well as the line $\Gamma _{d}$ where the
external sources are:%
\begin{equation}
\left. 
\begin{array}{c}
\Omega =\{\mathbf{x}:-A<x,y<A,a<z<b\}, \\ 
\partial _{1}\Omega =\left\{ \mathbf{x}:-A<x,y<A,z=a\right\} ,\text{ }%
\partial _{2}\Omega =\left\{ \mathbf{x}:-A<x,y<A,z=b\right\} , \\ 
\partial _{3}\Omega =\left\{ x=\pm A,y\in (-A.A),z\in \left( a,b\right)
\right\} \cup  \\ 
\cup \left\{ y=\pm A,x\in (-A,A),z\in \left( a,b\right) \right\} , \\ 
\Gamma _{d}=\{\mathbf{x}_{\alpha }=(\alpha ,0,0):\alpha \in \lbrack -d,d]\}.%
\end{array}%
\right.   \label{2.2}
\end{equation}%
Hence, $\Gamma _{d}$ is a part of the $x-$axis. By (\ref{2.1}) and (\ref{2.2}%
) $\Gamma _{d}\cap \overline{\Omega }=\varnothing $.

Let the points of external sources $\mathbf{x}_{\alpha }\in \Gamma _{d}$.
Let $\epsilon >0$ be a sufficiently small number. To avoid dealing with
singularities, we model the $\delta \left( \mathbf{x}\right) -$function as: 
\begin{equation}
f\left( \mathbf{x}\right) =C_{\epsilon }\left\{ 
\begin{array}{cc}
\exp \left( \frac{\left\vert \mathbf{x}\right\vert ^{2}}{\epsilon
^{2}-\left\vert \mathbf{x}\right\vert ^{2}}\right) , & \left\vert \mathbf{x}
\right\vert <\epsilon , \\ 
0, & \left\vert \mathbf{x}\right\vert \geq \epsilon ,%
\end{array}
\right.  \label{2.7}
\end{equation}%
where the constant $C_{\epsilon }$ is such that 
\begin{equation}
C_{\epsilon }\int_{\left\vert \mathbf{x}\right\vert <\epsilon }\exp \left( 
\frac{\left\vert \mathbf{x}\right\vert ^{2}}{\epsilon ^{2}-\left\vert 
\mathbf{x}\right\vert ^{2}}\right) d\mathbf{x}=1.  \label{2.8}
\end{equation}%
Hence, the function $f\left( \mathbf{x}-\mathbf{x}_{\alpha }\right) =f\left(
x-\alpha ,y,z\right) \in C^{\infty }\left( \mathbb{R}^{3}\right) $ plays the
role of the source function for the point source $\left\{ \mathbf{x}_{\alpha
}\right\} $. We choose $\epsilon $ so small that 
\begin{equation}
f\left( \mathbf{x}-\mathbf{x}_{\alpha }\right) =0,\quad \forall \mathbf{x}
\in \overline{\Omega },\quad \forall \mathbf{x}_{\alpha }\in \Gamma _{d}.
\label{2.9}
\end{equation}

Let $\Gamma (\mathbf{x},\mathbf{x}_{0})$ be the geodesic line generated by
the Riemannian metric $\sqrt{\varepsilon _{r}\left( \mathbf{x}\right) }%
\left\vert d\mathbf{x}\right\vert $ and connecting the source $\mathbf{x}%
_{0}\in \mathbb{R}^{3}$ with an arbitrary point $\mathbf{x}\in \mathbb{R}^{3}
$,

\begin{equation}
\Gamma (\mathbf{x},\mathbf{x}_{0})=\text{argmin}\left\{ 
\begin{array}{c}
\int\limits_{\gamma }\sqrt{\varepsilon _{r}\left( \mathbf{\xi }\left(
t\right) \right) }dt,\text{ where }\gamma \left( t\right) :[0,1]\rightarrow 
\mathbb{R}^{3} \\ 
\text{is a smooth map with }\gamma \left( 0\right) =\mathbf{x}_{0},\text{ }
\gamma \left( 0\right) =\mathbf{x}.%
\end{array}
\right\}  \label{2.10}
\end{equation}%
Here $\varepsilon _{r}\left( \mathbf{x}\right) $ is the spatially
distributed dielectric constant of the medium, $1/\sqrt{\varepsilon
_{r}\left( \mathbf{x}\right) }$ is the dimensionless speed of light. We
assume that the function $\varepsilon _{r}\left( \mathbf{x}\right) $
satisfies the following conditions: 
\begin{align}
\varepsilon _{r}(\mathbf{x})& \in C^{3}(\mathbb{R}^{3}),  \label{2.11} \\
\varepsilon _{r}(\mathbf{x})& =1,\quad \mathbf{x}\in \{\mathbf{x}\in \mathbb{%
\ R}^{3}\big\rvert\ \rvert x\rvert \geq A,\rvert y\rvert \geq A\}\cup \{ 
\mathbf{x}\in \mathbb{R}^{3}\big \rvert\ z\leq a\},  \label{2.12} \\
\partial _{z}\varepsilon _{r}(\mathbf{x})& \geq 0,\quad \mathbf{x}\in 
\mathbb{R}^{3}.  \label{2.13}
\end{align}%
Let $\tau (\mathbf{x},\mathbf{x}_{0})$ be the first time of arrival at the
point $\mathbf{x}$ of light generated at the point $\mathbf{x}_{0}.$ Then 
\cite[Chapter 3]{Rom}%
\begin{equation}
\tau (\mathbf{x},\mathbf{x}_{0})=\int\limits_{\Gamma (\mathbf{x},\mathbf{x}
_{0})}\sqrt{\varepsilon _{r}\left( \mathbf{\xi }\left( \sigma \right)
\right) }d\sigma ,  \label{2.130}
\end{equation}%
where $d\sigma $ is the element of the Euclidian length. For $\mathbf{x}\neq 
\mathbf{x}_{0}$ the function $\tau (\mathbf{x},\mathbf{x}_{0})$ is twice
continuously differentiable with respect to $\mathbf{x},\mathbf{x}_{0}$ and
is the solution of the eikonal equation \cite[Chapter 3]{Rom} 
\begin{equation}
\mid \nabla _{\mathbf{x}}\tau (\mathbf{x},\mathbf{x}_{0})\mid
^{2}=\varepsilon _{r}(\mathbf{x}),\text{ }\tau (\mathbf{x},\mathbf{x}
_{0})=O\left( \mid \mathbf{x}-\mathbf{x}_{0}\mid \right) ,\quad \mathbf{x}
\rightarrow \mathbf{x}_{0}.  \label{2.14}
\end{equation}

We assume everywhere below that the geodesic lines are regular \cite[Chapter
3]{Rom}:

\textbf{Regularity Assumption.}\emph{\ Any two points }$\mathbf{x,x}_{0}\in 
\mathbb{R}^{3}$\emph{\ can be connected by a single geodesic line }$\Gamma (%
\mathbf{x,x}_{0}).$\emph{\ }

A sufficient condition guaranteeing the regularity of geodesic lines can be
found in \cite{Rom2}. Let $\mu _{a}\left( \mathbf{x}\right) $ and $\mu _{s}(%
\mathbf{x})$ be the absorption and scattering coefficients of light
respectively and let%
\begin{equation}
\left. 
\begin{array}{c}
\mu _{a}\left( \mathbf{x}\right) ,\mu _{s}(\mathbf{x})\geq 0,\quad \mu
_{a}\left( \mathbf{x}\right) ,\mu _{s}(\mathbf{x})\in C^{1}\left( \mathbb{R}%
^{3}\right) , \\ 
\mu _{a}\left( \mathbf{x}\right) =\mu _{s}(\mathbf{x})=0,\quad \mathbf{x}\in 
\mathbb{R}^{3}\setminus \Omega , \\ 
a\left( \mathbf{x}\right) =\mu _{a}\left( \mathbf{x}\right) +\mu _{s}(%
\mathbf{x}).%
\end{array}%
\right.   \label{2.17}
\end{equation}%
The function $a\left( \mathbf{x}\right) $ is the attenuation coefficient. By
(\ref{2.17}) 
\begin{equation}
a\left( \mathbf{x}\right) \geq 0,\ \mathbf{x}\in \mathbb{R}^{3},\text{ }%
a\left( \mathbf{x}\right) \in C^{1}\left( \mathbb{R}^{3}\right) ,\quad
a\left( \mathbf{x}\right) =0,\ \mathbf{x}\in \mathbb{R}^{3}\setminus \Omega .
\label{2.20}
\end{equation}

Let $\widetilde{A}=\max (A,d)$. Introduce three domains $G,G_{a}^{+}$ and $%
G_{a}^{-},$ 
\begin{equation}
G=\left\{ \mathbf{x}:-\widetilde{A}<x,y<\widetilde{A},z\in (0,b)\right\}
,G_{a}^{+}=G\cup \{z>a\},G_{a}^{-}=G\setminus G_{a}^{+}.  \label{2.21}
\end{equation}%
Below we write sometimes $u(\mathbf{x},\alpha )$ instead of $u(\mathbf{x},%
\mathbf{x}_{\alpha })$.

\textbf{The Forward Problem.} \emph{Find the solution} $u(\mathbf{x},\alpha
)\in C^{1}\Big(G\times \left[ -d,d\right] \Big)$ \emph{of the following
problem:}%
\begin{equation}
\left. 
\begin{array}{c}
\left( \nabla _{\mathbf{x}}\tau (\mathbf{x},\mathbf{x}_{\alpha })/\sqrt{%
\varepsilon _{r}(\mathbf{x})}\right) \cdot \nabla _{\mathbf{x}}u(\mathbf{x}%
,\alpha )+a(\mathbf{x})u(\mathbf{x},\alpha )= \\ 
=\mu _{s}(\mathbf{x})\int_{\Gamma _{d}}K(\mathbf{x},\alpha ,\beta )u(\mathbf{%
\ x},\beta )d\beta +f(\mathbf{x}-\mathbf{x}_{\alpha }),\text{ }\mathbf{\ x}%
\in G,\mathbf{x}_{\alpha }\in \Gamma _{d},%
\end{array}%
\right.  \label{2.22}
\end{equation}%
\begin{equation}
u(\mathbf{x}_{\alpha },\mathbf{x}_{\alpha })=0\text{ for }\mathbf{x}_{\alpha
}\in \Gamma _{d}.  \label{2.23}
\end{equation}

\textbf{Definition 2.1.} \emph{We call equation (\ref{2.22}) the Riemannian
Radiative Transfer Equation (RRTE).}

In (\ref{2.22}), (\ref{2.23}) $u(\mathbf{x},\alpha )$ denotes the
steady-state radiance at the point $\mathbf{x}$ generated by the source
function $f\left( \mathbf{x}-\mathbf{x}_{\alpha }\right) $. The kernel $K(%
\mathbf{x},\alpha ,\beta )$ of the integral operator in (\ref{2.22}) is
called the \textquotedblleft phase function" \cite{Heino},%
\begin{equation}
\left. 
\begin{array}{c}
K(\mathbf{x},\alpha ,\beta )\geq 0,\quad \mathbf{x}\in \overline{\Omega }
;\quad \alpha ,\beta \in \left[ -d,d\right] , \\ 
K(\mathbf{x},\alpha ,\beta )\in C^{1}\left( \overline{\Omega }\times \left[
-d,d\right] ^{2}\right) .%
\end{array}
\right.  \label{2.24}
\end{equation}

\textbf{Coefficient Inverse Problem.} \emph{Let the function} $u\left( 
\mathbf{x},\alpha \right) $ $\in C^{1}(\overline{\Omega }$ $\times \left[
-d,d\right] )$ \emph{be the solution of the Forward Problem. Assume that the
coefficient} $a\left( \mathbf{x}\right) $ \emph{of equation \eqref{2.22} is
unknown. Determine the function} $a\left( \mathbf{x}\right) $, \emph{\
assuming that the following function} $g\left( \mathbf{x},\alpha \right) $ 
\emph{is known:} 
\begin{equation}
g\left( \mathbf{x},\alpha \right) =u\left( \mathbf{x},\alpha \right) ,\quad
\forall \mathbf{x}\in \partial \Omega \diagdown \partial _{1}\Omega ,\quad
\forall \alpha \in \left( -d,d\right) .  \label{2.42}
\end{equation}

\section{Existence and Uniqueness Theorem for the Forward Problem}

\label{sec:3}

Consider the unit tangent vector $\nu $ to the geodesic line $\Gamma (%
\mathbf{x},\mathbf{x}_{\alpha })$ at the point $\mathbf{x}$ \cite[Chapter 3]%
{Rom} 
\begin{equation*}
\nu =\nabla _{\mathbf{x}}\tau \left( \mathbf{x},\mathbf{x}_{\alpha }\right) /%
\sqrt{\varepsilon _{r}\left( \mathbf{x}\right) }.
\end{equation*}%
Hence, the directional derivative $D_{\nu }q$ of an appropriate function $q(%
\mathbf{x},\alpha )$ in the direction of the vector $\nu $ is 
\begin{equation}
D_{\nu }q=\frac{\nabla _{\mathbf{x}}\tau (\mathbf{x},\mathbf{x}_{\alpha })}{%
\sqrt{\varepsilon _{r}(\mathbf{x})}}\cdot \nabla _{\mathbf{x}}q(\mathbf{x}%
,\alpha ).  \label{2.15}
\end{equation}%
Hence, if the function $q(\mathbf{x},\mathbf{x}_{\alpha })$ solves problem (%
\ref{2.150}), then $q$ is given by formula (\ref{2.16}), where 
\begin{equation}
\frac{\nabla _{\mathbf{x}}\tau (\mathbf{x},\mathbf{x}_{\alpha })}{\sqrt{%
\varepsilon _{r}(\mathbf{x})}}\cdot \nabla _{\mathbf{x}}q(\mathbf{x},\mathbf{%
\ \ x}_{\alpha })=a(\mathbf{x}),\quad q(\mathbf{x}_{\alpha },\mathbf{x}%
_{\alpha })=0,  \label{2.150}
\end{equation}%
\begin{equation}
q(\mathbf{x},\mathbf{x}_{\alpha })=\int\limits_{\Gamma (\mathbf{x},\mathbf{x}%
_{\alpha })}a\left( \mathbf{\xi }(\sigma )\right) d\sigma .  \label{2.16}
\end{equation}%
Let 
\begin{equation}
p(\mathbf{x},\mathbf{x}_{\alpha })=\exp \left( \int\limits_{\Gamma (\mathbf{x%
},\mathbf{x}_{\alpha })}a\left( \mathbf{\xi }(\sigma )\right) d\sigma
\right) .  \label{2.27}
\end{equation}%
Then (\ref{2.15})-(\ref{2.27}) imply: 
\begin{equation}
D_{\nu }p=a\left( \mathbf{x}\right) p.  \label{2.28}
\end{equation}%
Multiply both sides of equation (\ref{2.22}) by $p$ and use (\ref{2.15})-(%
\ref{2.28}). Note that by (\ref{2.9}) and (\ref{2.17}) $p(\mathbf{x},\mathbf{%
x}_{\alpha })f(\mathbf{x}-\mathbf{x}_{\alpha })=f(\mathbf{\ \ x}-\mathbf{x}%
_{\alpha }).$ We obtain%
\begin{equation}
\left. 
\begin{array}{c}
pD_{\nu }u+a\left( \mathbf{x}\right) pu=\mu _{s}(\mathbf{x}%
)p\int\limits_{\Gamma _{d}}K(\mathbf{x},\alpha ,\beta )u(\mathbf{x},\beta
)d\beta +f(\mathbf{x}-\mathbf{x}_{\alpha }), \\ 
pD_{\nu }u+a\left( \mathbf{x}\right) pu=D_{\nu }\left( pu\right) -a\left( 
\mathbf{x}\right) pu+a\left( \mathbf{x}\right) pu=D_{\nu }\left( pu\right) ,
\\ 
D_{\nu }\left( pu\right) =\mu _{s}(\mathbf{x})p\int\limits_{\Gamma _{d}}K(%
\mathbf{x},\alpha ,\beta )u(\mathbf{x},\beta )d\beta +f(\mathbf{x}-\mathbf{x}%
_{\alpha }).%
\end{array}%
\right.   \label{2.29}
\end{equation}

Let the equation of the geodesic line $\Gamma (\mathbf{x},\mathbf{x}_{\alpha
})$ be $\mathbf{\xi =\xi }\left( \sigma ,\alpha \right) \in \Gamma (\mathbf{x%
},\mathbf{x}_{\alpha })$, where $\sigma $ is the Euclidean length of the
part $\Gamma _{\mathbf{\xi }}(\mathbf{x},\mathbf{x}_{\alpha })$ of the curve 
$\Gamma (\mathbf{x},\mathbf{x}_{\alpha }),$\ which connects points $\mathbf{%
\ \ \xi }$ and $\mathbf{x}_{\alpha }$. Integrating the last line of (\ref%
{2.29}) along the vector $\nu $ and taking into account the initial
condition (\ref{2.23}), we obtain for $\mathbf{x}\in G,\mathbf{x}_{\alpha
}\in \Gamma _{d}$%
\begin{equation}
\left. 
\begin{array}{c}
u(\mathbf{x},\mathbf{x}_{\alpha })=u_{0}(\mathbf{x},\mathbf{x}_{\alpha
})+p^{-1}(\mathbf{x},\mathbf{x}_{\alpha })\times \\ 
\times \int\limits_{\Gamma (\mathbf{x},\mathbf{x}_{\alpha })}p(\mathbf{\xi }%
(\sigma ,\alpha ),\mathbf{x}_{\alpha })\mu _{s}(\mathbf{\xi }(\sigma ,\alpha
))\left( \int\limits_{\Gamma _{d}}K(\mathbf{\xi }(\sigma ,\alpha ),\alpha
,\beta )u(\mathbf{\xi }(\sigma ,\alpha ),\beta )d\beta \right) d\sigma , \\ 
u_{0}(\mathbf{x},\mathbf{x}_{\alpha })=p^{-1}(\mathbf{x},\mathbf{x}_{\alpha
})\int\limits_{\Gamma \left( \mathbf{x},\mathbf{x}_{\alpha }\right) }f(%
\mathbf{\xi }(\sigma ,\alpha )-\mathbf{x}_{\alpha })d\sigma .%
\end{array}%
\right.  \label{2.32}
\end{equation}

Thus, we conclude that the solution of the Forward Problem (\ref{2.22}), (%
\ref{2.23}) is equivalent to the solution of integral equation (\ref{2.32}).

\textbf{Theorem 1.} \emph{Assume that conditions (\ref{2.17}) and (\ref{2.24}%
) hold. Then there exists unique solution }$u(\mathbf{x},\alpha )\in C^{1}%
\Big(G\times \left[ d,d\right] \Big)$\emph{\ of problem ( \ref{2.22}), (\ref%
{2.23}). Furthermore, the following inequality is valid: }%
\begin{align}
& u(\mathbf{x},\alpha )\geq m>0\text{ for }(\mathbf{x},\alpha )\in \left( 
\overline{G}_{a}^{+}\times \lbrack -d,d]\right) ,  \label{2.33} \\
& \hspace{1cm}m=\min_{(\mathbf{x},\alpha )\in \left( \overline{G}%
_{a}^{+}\times \lbrack -d,d]\right) }u_{0}(\mathbf{x},\alpha ),  \label{2.34}
\end{align}%
\emph{where the domain }$G_{a}^{+}$\emph{\ is defined in (\ref{2.21}).
Solution of problem (\ref{2.22}), (\ref{2.23}) is equivalent to the solution
of equation (\ref{2.32}). }

\textbf{Proof. }The equivalency was proven above in this section.\textbf{\ }%
Let $\mathbf{x}^{\ast }$ be the intersection point of the geodesic line $%
\Gamma (\mathbf{x},\mathbf{x}_{\alpha })$ with plane $\left\{ z=a\right\} $.
Note that by (\ref{2.1}), (\ref{2.2}) and (\ref{2.12}) $\Gamma (\mathbf{x}%
^{\ast },\mathbf{x}_{\alpha })$ is an interval of a straight line. Since by (%
\ref{2.1}), (\ref{2.17}) and (\ref{2.21}) $\mu _{s}(\mathbf{x})=0$ for $%
\mathbf{x}\in G_{a}^{-}$, then the first two lines of (\ref{2.32}) can be
rewritten as:%
\begin{equation}
\left. 
\begin{array}{c}
u(\mathbf{x},\mathbf{x}_{\alpha })=u_{0}(\mathbf{x},\mathbf{x}_{\alpha })+
\\ 
+p^{-1}(\mathbf{x},\mathbf{x}_{\alpha })\int\limits_{\Gamma (\mathbf{x},%
\mathbf{x}^{\star })}\left( \int\limits_{\Gamma _{d}}\widehat{K}(\mathbf{\xi 
}(\sigma ,\alpha ),\alpha ,\beta )u(\mathbf{\xi }(\sigma ,\alpha ),\beta
)d\beta \right) d\sigma ,%
\end{array}%
\right.  \label{2.35}
\end{equation}%
where the function $u_{0}(\mathbf{x},\mathbf{x}_{\alpha })$ is given in the
third line of (\ref{2.32}) and 
\begin{equation}
\widehat{K}(\mathbf{x},\alpha ,\beta )=p(\mathbf{x},\mathbf{x}_{\alpha })\mu
_{s}(\mathbf{x})K(\mathbf{x},\alpha ,\beta ).  \label{2.36}
\end{equation}

Consider now equations of the geodesic lines. Denote 
\begin{equation}
q_{1}=\tau _{x}(\mathbf{x},\mathbf{x}_{\alpha }),\text{ }q_{2}=\tau _{y}(%
\mathbf{x},\mathbf{x}_{\alpha }),\text{ }q_{3}=\tau _{z}(\mathbf{x},\mathbf{x%
}_{\alpha }).  \label{2.360}
\end{equation}%
Then formulas (3.4) and (3.7) of \cite[Chapter 3]{Rom} imply that equations
of geodesic lines are:%
\begin{equation}
\frac{dx}{ds}=\frac{q_{1}}{\varepsilon _{r}},\frac{dy}{ds}=\frac{q_{2}}{%
\varepsilon _{r}},\text{ }\frac{dz}{ds}=\frac{q_{3}}{\varepsilon _{r}},\text{
}  \label{2.361}
\end{equation}%
\begin{equation*}
\frac{dq_{1}}{ds}=\frac{\partial _{x}\varepsilon _{r}}{2\varepsilon _{r}},%
\text{ }\frac{dq_{2}}{ds}=\frac{\partial _{y}\varepsilon _{r}}{2\varepsilon
_{r}},\text{ }\frac{dq_{3}}{ds}=\frac{\partial _{z}\varepsilon _{r}}{%
2\varepsilon _{r}},
\end{equation*}%
where $ds=\sqrt{\varepsilon _{r}(\mathbf{x}(\sigma ))}d\sigma $ is the
element of the Riemannian length. In the integral (\ref{2.35}), 
\begin{equation}
\mathbf{x}(\sigma ,\alpha )=\left( x(\sigma ,\alpha ),y(\sigma ,\alpha
),z(\sigma ,\alpha )\right) \in \Omega .  \label{2.0036}
\end{equation}%
It follows from (\ref{2.13}) and \cite[Lemma 5.1]{KR} that there exists a
number $c>0$ such that 
\begin{equation}
\tau _{z}\left( \mathbf{x},\mathbf{x}_{\alpha }\right) \geq c.  \label{2.362}
\end{equation}%
Hence, a combination of equation (\ref{2.360}) with the last equation in (%
\ref{2.361}) implies: 
\begin{equation}
\partial _{s}z(s,\alpha )>0\text{ and \ }\partial _{\sigma }z(\sigma ,\alpha
)>0.  \label{2.036}
\end{equation}%
Consider the equation of the geodesic line $\Gamma (\mathbf{x},\mathbf{x}%
^{\star })$ in the form: 
\begin{equation}
\mathbf{\xi }(\sigma ,\alpha )=(\xi (\sigma ,\alpha ),\eta (\sigma ,\alpha
),\zeta (\sigma ,\alpha ))  \label{2.00036}
\end{equation}%
Change variables in (\ref{2.00036}) by replacing the variable $\sigma $ with
the variable $\zeta =\zeta (\sigma ,\alpha )$. Let $\sigma =\sigma (\zeta
,\alpha )$ be the inverse function. Then the equation of the geodesic line $%
\Gamma (\mathbf{x},\mathbf{x}^{\star })$ can be rewritten as 
\begin{equation*}
\mathbf{\xi }=\widehat{\mathbf{\xi }}(\zeta ,\alpha )=\mathbf{\xi }(\sigma
(\zeta ,\alpha ),\alpha )=(\xi (\sigma (\zeta ,\alpha ),\alpha ),\eta
(\sigma (\zeta ,\alpha ),\alpha ),\zeta ),\quad \zeta \in (a,z).
\end{equation*}%
By (\ref{2.0036}), (\ref{2.036}) and (\ref{2.00036}) the inverse function $%
\sigma =\sigma (\zeta ,\alpha )$ is monotonically increasing with respect to 
$\zeta $ along\textbf{\ }the geodesic line $\Gamma (x,x^{\star })$, i.e. $%
\partial _{\zeta }\sigma \left( \zeta ,\alpha \right) >0.$ Hence, we change
variables in the integral of (\ref{2.35}) as: $\sigma \Leftrightarrow \zeta
=\zeta (\sigma ,\alpha )$. Then equation (\ref{2.35}) can be rewritten as:%
\begin{equation}
u(\mathbf{x},\mathbf{x}_{\alpha })=u_{0}(\mathbf{x},\mathbf{x}_{\alpha
})+\int\limits_{a}^{z}\left( \int\limits_{\Gamma _{d}}\widetilde{K}(\mathbf{x%
},\widehat{\mathbf{\xi }}(\zeta ,\alpha ),\alpha ,\beta ,\zeta )u(\widehat{%
\mathbf{\xi }}(\zeta ,\alpha ),\beta )d\beta \right) d\zeta ,  \label{2.37}
\end{equation}%
where $\mathbf{x}\in G_{a}^{+},\mathbf{x}_{\alpha }\in \Gamma _{d}$ and by (%
\ref{2.36})%
\begin{equation}
\widetilde{K}(\mathbf{x},\mathbf{\xi },\alpha ,\beta ,\zeta )=\frac{1}{p(%
\mathbf{x},\mathbf{x}_{\alpha })}\widehat{K}(\mathbf{\xi },\alpha ,\beta
)\partial _{\zeta }\sigma (\zeta ,\alpha ).  \label{2.38}
\end{equation}%
Since we have the integral 
\begin{equation*}
\int\limits_{a}^{z}\left( ...\right) d\zeta 
\end{equation*}%
in equation (\ref{2.37}), then this is the integral equation of the Volterra
type. It follows from (\ref{2.11}), (\ref{2.24}), (\ref{2.27}), (\ref{2.36}%
)-(\ref{2.361}) and (\ref{2.38}) that the kernel of equation (\ref{2.37}) is
a non negative continuously differentiable function of its variables $\left( 
\mathbf{x},\alpha ,\beta ,z\right) \mathbf{\in }\overline{G}_{a}^{+}$ $%
\times \overline{\Gamma }_{d}$ $\times \overline{\Gamma }_{d}\times \lbrack
a,b]$. Hence, there exists a number $K_{0}>0$ such that in (\ref{2.37}), (%
\ref{2.38}) 
\begin{equation}
0\leq \widetilde{K}(\mathbf{x},\mathbf{\xi },\alpha ,\beta ,\zeta )\leq
K_{0}<\infty \text{ in (\ref{2.37}).}  \label{2.39}
\end{equation}

Since equation (\ref{2.37}) is of the Volterra type, then its solution can
be obtained iteratively as:%
\begin{equation}
\left. 
\begin{array}{c}
u_{n}(\mathbf{x},\mathbf{x}_{\alpha
})=\int\limits_{a}^{z}\int\limits_{\Gamma _{d}}\widetilde{K}(\mathbf{x},%
\widehat{\mathbf{\xi }}(\zeta ,\alpha ),\alpha ,\beta ,\zeta )u_{n-1}(%
\widehat{\mathbf{\xi }}(\zeta ,\alpha ),\beta )d\beta d\zeta , \\ 
u(\mathbf{x},\mathbf{x}_{\alpha })=\sum_{n=0}^{\infty }u_{n}(\mathbf{x},%
\mathbf{x}_{\alpha })\text{.}%
\end{array}%
\right.   \label{2.40}
\end{equation}%
It follows from (\ref{2.7}), (\ref{2.8}), (\ref{2.34}) and (\ref{2.37})-(\ref%
{2.40}) that%
\begin{equation}
\left. 
\begin{array}{c}
m\leq u(\mathbf{x},\mathbf{x}_{\alpha })\leq \left[ \max_{(\mathbf{x},\alpha
)\in \left( G_{a}^{+}\times \lbrack -d,d]\right) \ }u_{0}\left( \mathbf{x},%
\mathbf{x}_{\alpha }\right) \right] \times  \\ 
\times \dsum\limits_{n=0}^{\infty }\left( 2dK_{0}(z-a)\right) ^{n}/n!,\text{ 
}\mathbf{x}\in G_{a}^{+},%
\end{array}%
\right.   \label{2.41}
\end{equation}%
where numbers $m$ and $K_{0}$ are defined in (\ref{2.34}) and (\ref{2.39})
respectively. Estimate (\ref{2.33}) follows from (\ref{2.41}). Obviously the
series of first derivatives of terms of (\ref{2.40}) with respect to any of
variables $x,y,z,\alpha $ also converges absolutely. Hence the function $u(%
\mathbf{x},\mathbf{x}_{\alpha })$ in (\ref{2.40}) belongs to $C^{1}\Big(%
\overline{G}_{a}^{+}\times \overline{\Gamma }_{d}\Big).$ We set%
\begin{equation*}
u(\mathbf{x},\mathbf{x}_{\alpha })=\left\{ 
\begin{array}{c}
\text{the right hand side of (\ref{2.40}) for }(\mathbf{x},\mathbf{x}%
_{\alpha })\in G_{a}^{+}\times \Gamma _{d}, \\ 
u_{0}(\mathbf{x},\mathbf{x}_{\alpha })\text{ for }(\mathbf{x},\mathbf{x}%
_{\alpha })\in G_{a}^{-}\times \Gamma _{d}.%
\end{array}%
\right. 
\end{equation*}

Hence, the so defined function $u(\mathbf{x},\mathbf{x}_{\alpha })\in
C^{1}\left( \overline{G}\times \overline{\Gamma }_{d}\right) .$ Thus, we
have proven the existence of the solution $u(\mathbf{x},\alpha )\in
C^{1}\left( \overline{G}\times \left[ d,d\right] \right) $ of the Forward
Problem (\ref{2.22}), (\ref{2.23}) as well as estimate (\ref{2.33}). To
prove uniqueness, one should set in (\ref{2.37}) $u_{0}\left( \mathbf{x},%
\mathbf{x}_{\alpha }\right) \equiv 0$ and then proceed in the classical way
of the proof of the uniqueness of the Volterra integral equation of the
second kind. $\square $

\textbf{Remark 3.1. }\emph{It follows from (\ref{2.32}) and Theorem 1 that
one can solve Forward Problem via the solution of the linear integral
equation in (\ref{2.32}). This is how we solve the forward problem (\ref%
{2.22}), (\ref{2.23}) in the numerical section 6 to generate the data for
the inverse problem. }

\section{Convexification Numerical Method for the Coefficient Inverse Problem%
}

\label{sec:4}

\subsection{An integral differential equation without the unknown
coefficient $a\left( \mathbf{x}\right) $}

\label{sec:4.1}

By (\ref{2.7})-(\ref{2.9}) equation (\ref{2.22}) can be rewritten as:%
\begin{equation}
\left. 
\begin{array}{c}
\left( \nabla _{\mathbf{x}}\tau (\mathbf{x},\mathbf{x}_{\alpha })/\sqrt{%
\varepsilon _{r}(\mathbf{x})}\right) \cdot \nabla _{\mathbf{x}}u(\mathbf{x}%
,\alpha )+a(\mathbf{x})u(\mathbf{x},\alpha )= \\ 
=\mu _{s}(\mathbf{x})\int\limits_{\Gamma _{d}}K(\mathbf{x},\alpha ,\beta )u(%
\mathbf{x},\beta )d\beta ,\text{ }\left( \mathbf{x},\alpha \right) \in
\Omega \times \left( -d,d\right) .%
\end{array}%
\right.  \label{3.2}
\end{equation}%
It follows from (\ref{2.33}) that we can consider a new function $v(\mathbf{x%
},\alpha ),$ 
\begin{equation}
v(\mathbf{x},\alpha )=\ln u(\mathbf{x},\alpha ),\quad \left( \mathbf{x}%
,\alpha \right) \in \Omega \times \left( -d,d\right) .  \label{3.1}
\end{equation}%
By (\ref{3.1}) $u(\mathbf{x},\alpha )=e^{v(\mathbf{x},\alpha )}.$
Substituting this in (\ref{3.2}), we obtain for $\left( \mathbf{x},\alpha
\right) \in \Omega \times \left( -d,d\right) :$ 
\begin{equation}
\left( \nabla _{\mathbf{x}}\tau (\mathbf{x},\mathbf{x}_{\alpha })/\sqrt{%
\varepsilon _{r}(\mathbf{x})}\right) \cdot \nabla _{\mathbf{x}}v(\mathbf{x}%
,\alpha )+a(\mathbf{x})=e^{-v(\mathbf{x},\alpha )}\mu _{s}(\mathbf{x}%
)\int\limits_{\Gamma _{d}}K(\mathbf{x},\alpha ,\beta )e^{v(\mathbf{x},\beta
)}d\beta .  \label{3.4}
\end{equation}%
In particular, (\ref{3.4}) implies that we can calculate the function $a(%
\mathbf{x})$ by the following formula:%
\begin{equation}
\left. 
\begin{array}{c}
a(\mathbf{x})=-\dint\limits_{\Gamma _{d}}\left( \nabla _{\mathbf{x}}\tau (%
\mathbf{x},\mathbf{x}_{\alpha })/\sqrt{\varepsilon _{r}(\mathbf{x})}\right)
\cdot \nabla _{\mathbf{x}}v(\mathbf{x},\alpha )d\alpha + \\ 
+\dint\limits_{\Gamma _{d}}\left( e^{-v(\mathbf{x},\alpha )}\mu _{s}(\mathbf{%
x})\int\limits_{\Gamma _{d}}K(\mathbf{x},\alpha ,\beta )e^{v(\mathbf{x}%
,\beta )}d\beta \right) d\alpha .%
\end{array}%
\right.  \label{3.5}
\end{equation}

Hence, we now focus on the problem of the reconstruction of the function $v(%
\mathbf{x},\alpha )$ from the function $g(\mathbf{x},\alpha )$ given in (\ref%
{2.42}). We have 
\begin{equation}
\frac{\tau _{z}\left( \mathbf{x},\alpha \right) }{\sqrt{\varepsilon _{r}(%
\mathbf{x})}}v_{z}\left( \mathbf{x},\alpha \right) =\frac{\partial }{%
\partial z}\left( \frac{\tau _{z}}{\sqrt{\varepsilon _{r}}}v\right) -\frac{%
\partial }{\partial z}\left( \frac{\tau _{z}}{\sqrt{\varepsilon _{r}}}%
\right) v.  \label{3.7}
\end{equation}%
Introduce a new function $w\left( \mathbf{x},\alpha \right) $ and express $%
v\left( \mathbf{x},\alpha \right) $ through $v\left( \mathbf{x},\alpha
\right) ,$%
\begin{equation}
\left. 
\begin{array}{c}
w\left( \mathbf{x},\alpha \right) =\left( \tau _{z}\left( \mathbf{x},\alpha
\right) /\sqrt{\varepsilon _{r}\left( \mathbf{x}\right) }\right) v\left( 
\mathbf{x},\alpha \right) , \\ 
v\left( \mathbf{x},\alpha \right) =\left( \sqrt{\varepsilon _{r}\left( 
\mathbf{x}\right) }/\tau _{z}\left( \mathbf{x},\alpha \right) \right) .%
\end{array}%
\right.  \label{3.8}
\end{equation}%
It follows from (\ref{2.362}) that the second line of formula (\ref{3.8})
makes sense. Thus, (\ref{3.7}) becomes%
\begin{equation}
\frac{\tau _{z}}{\sqrt{\varepsilon _{r}}}v_{z}=w_{z}-\left[ \frac{\partial }{%
\partial z}\left( \frac{\tau _{z}}{\sqrt{\varepsilon _{r}}}\right) \frac{%
\sqrt{\varepsilon _{r}}}{\tau _{z}}\right] w.  \label{3.10}
\end{equation}%
Using (\ref{3.8}), transform other terms of the differential operator in (%
\ref{3.4}),%
\begin{equation}
\frac{\tau _{x}}{\sqrt{\varepsilon _{r}}}v_{x}=\frac{\tau _{x}}{\sqrt{%
\varepsilon _{r}}}\frac{\partial }{\partial x}\left( \frac{\sqrt{\varepsilon
_{r}}}{\tau _{z}}w\right) =\frac{\tau _{x}}{\tau _{z}}w_{x}+\left[ \frac{%
\tau _{x}}{\sqrt{\varepsilon _{r}}}\frac{\partial }{\partial x}\left( \frac{%
\sqrt{\varepsilon _{r}}}{\tau _{z}}\right) \right] w.  \label{3.11}
\end{equation}%
And similarly for $\left( \tau _{y}/\sqrt{\varepsilon _{r}}\right) v_{y}.$
Hence, (\ref{3.4}) becomes%
\begin{equation}
\left. 
\begin{array}{c}
w_{z}+\left( \tau _{x}w_{x}+\tau _{y}w_{y}\right) /\tau _{z}+ \\ 
+\left[ \left( \tau _{x}/\sqrt{\varepsilon _{r}}\right) \partial _{x}\left( 
\sqrt{\varepsilon _{r}}/\tau _{z}\right) +\left( \tau _{y}/\sqrt{\varepsilon
_{r}}\right) \partial _{y}\left( \sqrt{\varepsilon _{r}}/\tau _{z}\right)
-\left( \sqrt{\varepsilon _{r}}/\tau _{z}\right) \partial _{z}\left( \tau
_{z}/\sqrt{\varepsilon _{r}}\right) \right] w- \\ 
-\exp \left( -w\sqrt{\varepsilon _{r}}/\tau _{z}\right) (\mathbf{x},\alpha
)\mu _{s}(\mathbf{x})\int\limits_{\Gamma _{d}}K(\mathbf{x},\alpha ,\beta
)\exp \left( w\sqrt{\varepsilon _{r}}/\tau _{z}\right) (\mathbf{x},\beta
)d\beta = \\ 
=-a\left( \mathbf{x}\right) ,\text{ }\left( \mathbf{x},\alpha \right) \in
\Omega \times \left( -d,d\right) .%
\end{array}%
\right.  \label{3.12}
\end{equation}

Differentiate both sides of (\ref{3.12}) with respect to $\alpha $ and use $%
\partial _{\alpha }a(\mathbf{x})\equiv 0$. We obtain for $\left( \mathbf{x}%
,\alpha \right) \in \Omega \times \left( -d,d\right) :$%
\begin{equation}
\left. 
\begin{array}{c}
\partial _{\alpha }w_{z}+\partial _{\alpha }\left( \left( \tau
_{x}w_{x}+\tau _{y}w_{y}\right) /\tau _{z}\right) + \\ 
+\partial _{\alpha }\left\{ \left[ \left( \tau _{x}/\sqrt{\varepsilon _{r}}%
\right) \partial _{x}\left( \sqrt{\varepsilon _{r}}/\tau _{z}\right) +\left(
\tau _{y}/\sqrt{\varepsilon _{r}}\right) \partial _{y}\left( \sqrt{%
\varepsilon _{r}}/\tau _{z}\right) \right] w\right\} - \\ 
-\partial _{\alpha }\left[ \left( \sqrt{\varepsilon _{r}}/\tau _{z}\right)
\partial _{z}\left( \tau _{z}/\sqrt{\varepsilon _{r}}\right) w\right] - \\ 
-\partial _{\alpha }\left[ \exp \left( -w\sqrt{\varepsilon _{r}}/\tau
_{z}\right) (\mathbf{x},\alpha )\mu _{s}(\mathbf{x})\int\limits_{\Gamma
_{d}}K(\mathbf{x},\alpha ,\beta )\exp \left( w\sqrt{\varepsilon _{r}}/\tau
_{z}\right) (\mathbf{x},\beta )d\beta \right] =0.%
\end{array}%
\right.  \label{3.13}
\end{equation}%
The Dirichlet boundary condition for the function $w(\mathbf{x},\alpha )$
is: 
\begin{align}
& w(\mathbf{x},\alpha )=\frac{\tau _{z}(\mathbf{x},\alpha )}{\sqrt{%
\varepsilon _{r}\left( \mathbf{x}\right) }}\ln g_{1}(\mathbf{x},\alpha
),\quad (\mathbf{x},\alpha )\in \partial \Omega \times (-d,d),  \label{3.14}
\\
& \hspace{0cm}g_{1}(\mathbf{x},\alpha )=\left\{ 
\begin{array}{ll}
g(\mathbf{x},\alpha ), & \mathbf{x}\in \partial \Omega \diagdown \partial
_{1}\Omega ,\quad \alpha \in (-d,d), \\ 
u_{0}(\mathbf{x},\alpha ), & \mathbf{x}\in \partial _{1}\Omega ,\quad \alpha
\in (-d,d).%
\end{array}%
\right.  \label{3.15}
\end{align}

Thus, we develop below a numerical method to obtain an approximate solution $%
w(\mathbf{x},\alpha )$ of problem (\ref{3.13})-(\ref{3.15}).

\subsection{A special orthonormal basis in $L_{2}(-d,d)$}

\label{sec:4.2}

First, we introduce a special orthonormal basis in $L_{2}(-d,d),$ which was
first discovered in \cite{Klib2017}, also, see \cite[section 6.2.3]{KL}.
Consider a linearly independent set of functions $\{\alpha ^{n}e^{\alpha
}\}_{n=0}^{\infty }\subset L_{2}(-d,d)$, which is complete in $L_{2}(-d,d)$.
The Gram-Schmidt orthonormalization procedure being applied to this set,
results in the orthonormal basis $\{\Psi _{n}\left( \alpha \right)
\}_{n=0}^{\infty }$ in $L_{2}(-d,d)$. The Gram-Schmidt procedure is unstable
when it is applied to an infinite number of functions. However, we have not
seen an instability when applying it to a relatively small number of
functions for $n\in \left[ 0,12\right] .$ The same was observed in a number
of previous publications of this research group, see, e.g. \cite{Khoa,KTR}, 
\cite[Chapters 7,10,12]{KL}.

Let $\left[ ,\right] $ be the scalar product in $L_{2}(-d,d)$. Denote $%
b_{s,k}=\left[ Q_{s}^{\prime },Q_{k}\right] .$ Then \cite{Klib2017}, \cite[%
section 6.2.3]{KL}%
\begin{equation}
b_{s,k}=\left\{ 
\begin{array}{c}
1,s=k, \\ 
0,s>k.%
\end{array}
\right.  \label{3.150}
\end{equation}%
Consider the $N\times N$ matrix $B_{N}=\left( b_{s,k}\right) _{\left(
s,k\right) =\left( 0,0\right) }^{\left( N-1,N-1\right) }$. Then (\ref{3.150}%
) implies that $\det B_{N}=1,$ which means that this matrix is invertible.
In fact, the existence of the matrix $B_{N}^{-1}$ for each $N\geq 1$ is the
key property why the basis $\{Q_{n}\left( \alpha \right) \}_{n=0}^{\infty }$
was originally constructed in \cite{Klib2017}. Indeed, consider, for example
either the basis of standard orthonormal polynomials or the basis of
trigonometric functions. In each of these, the first function is an
identical constant, which means that the first raw of an analog of the
matrix $B_{N}$ is zero.

\subsection{A boundary value problem for a system of nonlinear PDEs}

\label{sec:4.3}

We assume that the functions $w(\mathbf{x},\alpha )$, $w_{\alpha }(\mathbf{x}%
,\alpha )$ can be represented as truncated Fourier-like series 
\begin{equation}
w(\mathbf{x},\alpha )=\sum\limits_{n=0}^{N-1}w_{n}(\mathbf{x})Q_{n}(\alpha
),\quad w_{\alpha }(\mathbf{x},\alpha )=\sum\limits_{n=0}^{N-1}w_{n}(\mathbf{%
\ \ x})Q_{n}^{\prime }(\alpha )  \label{3.16}
\end{equation}%
with unknown coefficients $\left\{ w_{n}(\mathbf{x})\right\} _{n=0}^{N-1}.$
Thus, we focus below on the computation of the $N-$D vector function 
\begin{equation}
V(\mathbf{x})=\left( w_{0},w_{1},\cdots ,w_{N-1}\right) ^{T}(\mathbf{x}).
\label{3.016}
\end{equation}

\textbf{Remarks 4.1}:

\begin{enumerate}
\item \emph{The representations (\ref{3.16}) mean that this is a version of
the Galerkin method. However, unlike classical well-posed forward problems
for PDEs, where Galerkin method is used and its convergence at }$%
N\rightarrow \infty $ \emph{is usually proven, we cannot prove convergence
of our inversion numerical procedure described below for }$N\rightarrow
\infty $\emph{. This is basically because of the ill-posed nature of our
CIP. Thus, we actually work below within the framework of an approximate
mathematical model. Then, however, the question can be raised whether this
model really works numerically. The answer is positive, and this answer is
obtained computationally in section 6. We observe that very similar
truncated series were used in some other above cited works on the
convexification, such as, e.g. \cite{Khoa,KTR}, \cite[Chapters 7,10]{KL},
and all of them have demonstrated good numerical performances. Likewise,
truncated Fourier series were used in works of other authors about CIPs,
such as, e.g. \cite{GN,Kab1,Kab2,Nov} and also without proofs of convergence
of inversion procedures at }$N\rightarrow \infty .$\emph{\ Those proofs were
not provided for the same reason as the one here: the ill-posed nature of
CIPs. }

\item \emph{Finally, we refer to subsection 3.4 of \cite{KTR} for more
arguments in support of those of item 1. In particular, these arguments
include the well known fact that the Huygens-Fresnel theory of the
diffraction in optics is not yet rigorously derived from the Maxwell's
equations, see, e.g. a classic textbook \cite[ pages 412, 413]{BW}.
Philosophically, this fact is similar with the discussion of item 1.}
\end{enumerate}

Substitute (\ref{3.16}) in (\ref{3.13}). Next, sequentially multiply the
obtained equation by $Q_{n}(\alpha )$, $n=0,...,N-1$ and integrate with
respect to $\alpha \in (-d,d).$ We obtain the following system of coupled
quasilinear integral differential equations 
\begin{equation}
B_{N}V_{z}(\mathbf{x})+A_{1}(\mathbf{x})V_{x}(\mathbf{x})+A_{2}(\mathbf{x}%
)V_{y}(\mathbf{x})+F\left( V(\mathbf{x}),\mathbf{x}\right) =0,\quad \mathbf{%
\ \ x}\in \Omega ,  \label{3.17}
\end{equation}%
where $A_{1}(\mathbf{x})$ and $A_{2}(\mathbf{x})$ are $N\times N$ matrices
and $F\left( V(\mathbf{x}),\mathbf{x}\right) $ is a certain vector function,
which depends nonlinearly on $V(\mathbf{x}).$ Explicit formulas for $A_{1}(%
\mathbf{x})$, $A_{2}(\mathbf{x})$ and $F\left( V(\mathbf{x}),\mathbf{x}%
\right) $ can be easily written. However, we do not present them here for
brevity. In addition, the boundary condition for the vector function $V(%
\mathbf{x})$ is: 
\begin{align}
& \hspace{1cm}V\left( \mathbf{x}\right) \mid _{\partial \Omega }=P(\mathbf{x}%
)=\left( p_{0},p_{1},\cdots ,p_{N-1}\right) ^{T}(\mathbf{x}),  \label{3.21}
\\
& p_{n}(\mathbf{x})=\int\limits_{-d}^{d}\left[ \frac{\tau _{z}(\mathbf{x}%
,\alpha )}{\sqrt{\varepsilon _{r}\left( \mathbf{x}\right) }}\ln \left[ g_{1}(%
\mathbf{x},\alpha )\right] \right] Q_{n}(\alpha )d\alpha ,\quad n=0,1,\cdots
,N-1.  \label{3.22}
\end{align}%
Thus, we now have to solve the boundary value problem (\ref{3.17})-(\ref%
{3.22}).

To numerically calculate the derivatives of $\nabla _{\mathbf{x}}\tau (%
\mathbf{x},\alpha )$ with respect to $\alpha $, we represent $\nabla _{%
\mathbf{x}}\tau (\mathbf{x},\alpha )$ via the truncated Fourier series with
respect to the above basis $\left\{ Q_{n}(\alpha )\right\} _{n=0}^{N-1}$ as: 
\begin{equation}
\nabla _{\mathbf{x}}\tau (\mathbf{x},\alpha )=\sum\limits_{n=0}^{N-1}\left(
\nabla _{\mathbf{x}}\tau \right) _{n}(\mathbf{x})Q_{n}(\alpha ).
\label{3.1001}
\end{equation}%
Then we use explicit formulas for functions $Q_{n}(\alpha )$ to get 
\begin{equation}
\partial _{\alpha }\left( \nabla _{\mathbf{x}}\tau \right)
=\sum\limits_{n=0}^{N-1}\left( \nabla _{\mathbf{x}}\tau \right) _{s}(\mathbf{%
x})Q_{n}^{\prime }(\alpha ).  \label{3.1002}
\end{equation}%
Then equations (\ref{3.1001}) and (\ref{3.1002}) are used in (\ref{3.17})-(%
\ref{3.22}). Thus, it follows from (\ref{3.13}) and (\ref{3.16})-(\ref%
{3.1002}) that%
\begin{equation}
\left\{ 
\begin{array}{c}
A_{1}(\mathbf{x}),A_{2}(\mathbf{x})\in C_{N^{2}}\left( \overline{\Omega }%
\right) ,\text{ and the vector function} \\ 
F\left( V\left( \mathbf{x}\right) ,\mathbf{x}\right) \text{ is continuously
differentiable} \\ 
\text{with respect to its arguments for }\mathbf{x}\in \overline{\Omega }.%
\end{array}%
\right.  \label{3.1003}
\end{equation}%
Here and below for any integer $k\geq 2$ and for any Banach space $B$ we
denote $B_{k}=B^{k}$ with the norm $\left\Vert f\right\Vert
_{B_{k}}^{2}=\left\Vert f_{1}\right\Vert _{B}^{2}+...+\left\Vert
f_{k}\right\Vert _{B}^{2},$ $\forall f=\left( f_{1},...,f_{k}\right) ^{T}\in
B_{k}.$

\subsection{Minimization problem}

\label{sec:4.4}

Let $R>0$ be an arbitrary number and the vector function $P\left( x\right) $
be the boundary condition in (\ref{3.21}). Define the set $S\left(
R,P\right) \subset H_{N}^{1}\left( \Omega \right) $ as:%
\begin{equation}
S\left( R,P\right) =\left\{ V\in H_{N}^{1}\left( \Omega \right) :V(\mathbf{x}
)\mid _{\partial \Omega }=P(\mathbf{x}),\left\Vert W\right\Vert
_{H_{N}^{1}\left( \Omega \right) }<R\right\} .  \label{3.23}
\end{equation}

To solve problem (\ref{3.17})-(\ref{3.22}), we solve the following
minimization problem:

\textbf{Minimization Problem 1.} \emph{Let }$\lambda \geq 1$\emph{\ be a
parameter. Minimize} \emph{the following weighted cost functional }$%
J_{\lambda }\left( V\right) $ \emph{\ on the set }$\overline{S\left(
R,P\right) }:$

\begin{equation}
\left. J_{\lambda }\left( V\right) =\left\Vert \left( B_{N}V_{z}+A_{1}( 
\mathbf{x})V_{x}(\mathbf{x})+A_{2}(\mathbf{x})V_{y}(\mathbf{x})+F\left( V( 
\mathbf{x}),\mathbf{x}\right) \right) e^{\lambda z}\right\Vert
_{L_{N}^{2}\left( \Omega \right) }^{2}.\right.  \label{3.240}
\end{equation}

\section{Convergence Analysis}

\label{sec:5}

We carry out the convergence analysis for a modified Minimization Problem 1.
To obtain this modification, we rewrite the differential operator in
functional (\ref{3.240}) via finite differences with respect to the
variables $x,y$ while leaving the conventional derivative with respect to $z$%
. We call this \textquotedblleft partial finite differences".

\subsection{Partial finite differences}

\label{sec:5.1}

Let $m>1$ be an integer. Let $A>0$ be the number in (\ref{2.2}). Consider
two partitions of the interval $\left( -A,A\right) $,%
\begin{equation}
\left. 
\begin{array}{c}
-A=x_{0}<x_{1}<\cdots <x_{m}=A,\quad x_{j+1}-x_{j}=h,\quad j=0,\cdots ,m-1,
\\ 
-A=y_{0}<y_{1}<\cdots <y_{m}=A,\quad y_{j+1}-y_{j}=h,\quad j=0,\cdots ,m-1.%
\end{array}%
\right.  \label{3.25}
\end{equation}%
We assume that 
\begin{equation}
h\geq h_{0}=const.>0.  \label{3.26}
\end{equation}%
Define the semidiscrete subset $\Omega ^{h}$ of the domain $\Omega $ as: 
\begin{align}
& \hspace{3cm}\Omega _{1}^{h}=\left\{ \left( x_{i},y_{j}\right) \right\}
_{i,j=0}^{m},  \label{3.27} \\
& \Omega ^{h}=\Omega _{1}^{h}\times \left( a,b\right) =\left\{
(x_{i},y_{j}):(x_{i},y_{j})\in \Omega _{1}^{h},z\in (a,b)\right\} .
\label{3.28}
\end{align}%
Below points $(x_{i},y_{j},z)\in \Omega ^{h}$ are denoted as $\mathbf{x}^{h}$%
. By (\ref{2.2}), (\ref{3.27}) and (\ref{3.28}) the boundary $\partial
\Omega ^{h}$ of the domain $\Omega ^{h}$ is: 
\begin{align*}
& \hspace{1.5cm}\partial \Omega ^{h}=\partial _{1}\Omega ^{h}\cup \partial
_{2}\Omega ^{h}\cup \partial _{3}\Omega ^{h}, \\
& \hspace{0.2cm}\partial _{1}\Omega ^{h}=\Omega _{1}^{h}\times \left\{
z=a\right\} ,\ \partial _{2}\Omega ^{h}=\Omega _{1}^{h}\times \left\{
z=b\right\} , \\
& \partial _{3}\Omega ^{h}=\left\{ \left( x_{0},y_{j},z\right) ,\left(
x_{m},y_{j},z\right) :z\in (a,b)\right\} .
\end{align*}%
Let the vector function $Y(\mathbf{x})\in C_{N}^{1}(\overline{\Omega })$.
Denote 
\begin{equation*}
Y^{h}(\mathbf{x}^{h})=Y(x_{i},y_{j},z),\quad \mathbf{x}^{h}=(x_{i},y_{j},z)%
\in \Omega ^{h}.
\end{equation*}%
Thus, $Y^{h}(\mathbf{x}^{h})$ is an $N-D$ vector function of discrete
variables $(x_{i},y_{j})\in \Omega _{1}^{h}$ and continuous variable $z\in
(a,b)$. Note that by (\ref{3.25}) the boundary terms at $\partial _{3}\Omega
^{h}$ of this vector function, which correspond to $Y(\mathbf{x})\mid
_{\partial _{3}\Omega ^{h}}$, are: 
\begin{equation*}
\left\{ Y(x_{0},y_{j},z)\right\} \cup \left\{ Y(x_{m},y_{j},z)\right\} \cup
\left\{ Y(x_{i},y_{0},z)\right\} \cup \left\{ Y(x_{i},y_{m},z)\right\}
,i,j=0,\cdots ,m.
\end{equation*}%
For two vector functions $Y^{\left( 1\right) }(\mathbf{x})=\left(
Y_{0}^{\left( 1\right) }(\mathbf{x}),\cdots ,Y_{N-1}^{\left( 1\right) }(%
\mathbf{x})\right) ^{T}$ and $Y^{\left( 2\right) }(\mathbf{x}%
)=(Y_{0}^{\left( 2\right) }(\mathbf{x})$, $\cdots $, $Y_{N-1}^{\left(
2\right) }(\mathbf{x}))^{T}$ their scalar product $Y^{\left( 1\right) }(%
\mathbf{x})\cdot Y^{\left( 2\right) }(\mathbf{x})$ is defined as the scalar
product in $\mathbb{R}^{N},$ and $\left( Y(\mathbf{x})\right) ^{2}=Y(\mathbf{%
x})\cdot Y(\mathbf{x})$. Respectively,%
\begin{equation}
\left. 
\begin{array}{c}
Y^{\left( 1\right) h}(\mathbf{x}^{h})\cdot Y^{\left( 2\right) h}(\mathbf{x}%
^{h})= \\ 
=\sum\limits_{n=0}^{N-1}\sum\limits_{\left( i,j\right) =\left( 1,1\right)
}^{\left( i,j\right) =\left( m-1,m-1\right) }Y_{n}^{\left( 1\right)
h}(x_{i},y_{j},z)Y_{n}^{\left( 2\right) h}(x_{i},y_{j},z), \\ 
\left( Y^{h}(\mathbf{x}^{h})\right) ^{2}=Y^{h}(\mathbf{x}^{h})\cdot Y^{h}(%
\mathbf{x}^{h}),\text{ }\left\vert Y^{h}(\mathbf{x}^{h})\right\vert =\sqrt{%
Y^{h}(\mathbf{x}^{h})\cdot Y^{h}(\mathbf{x}^{h})}.%
\end{array}%
\right.  \label{3.29}
\end{equation}%
We will use formulas (\ref{3.29}) everywhere below without further
mentioning. We exclude here boundary terms with $i,j=0$ and $i,j=m$ since we
work below with finite difference derivatives as defined in the next
paragraph.

We define finite difference derivatives of the semidiscrete $N-$D vector
function $Y^{h}(\mathbf{x}^{h})$ with respect to $x,y$ only at interior
points of the domain $\Omega ^{h}$ with $i,j=1,...,m-1$,%
\begin{equation}
\left. 
\begin{array}{c}
\partial _{x}Y^{h}\left( x_{i},y_{j},z\right) =Y^{h}\left(
x_{i},y_{j},z\right) _{x}=\left( Y^{h}\left( x_{i+1},y_{j},z\right)
-Y^{h}\left( x_{i-1},y_{j},z\right) \right) /\left( 2h\right) , \\ 
\partial _{x}Y^{h}\left( x_{i},y_{j},z\right) =Y^{h}\left(
x_{i},y_{j},z\right) _{x}=\left( Y^{h}\left( x_{i+1},y_{j},z\right)
-Y^{h}\left( x_{i-1},y_{j},z\right) \right) /\left( 2h\right) , \\ 
\partial _{y}Y^{h}\left( x_{i},y_{j},z\right) =Y^{h}\left(
x_{i},y_{j},z\right) _{y}=\left( Y^{h}\left( x_{i},y_{j+1},z\right)
-Y^{h}\left( x_{i},y_{j-1},z\right) \right) /\left( 2h\right) , \\ 
Y_{x}^{h}\left( \mathbf{x}^{h}\right) =\left\{ Y^{h}\left(
x_{i},y_{j},z\right) _{x}\right\} _{i,j=1}^{m-1},\text{ }Y_{y}^{h}\left( 
\mathbf{x}^{h}\right) =\left\{ Y^{h}\left( x_{i},y_{j},z\right) _{y}\right\}
_{i,j=1}^{m-1}.%
\end{array}
\right.  \label{3.31}
\end{equation}

We need semidiscrete analogs of spaces $C_{N^{2}}\left( \overline{\Omega }%
\right) ,H_{N}^{1}\left( \Omega \right) ,L_{N}^{2}\left( \Omega \right) $.
All three are defined using the same principle. Hence, we provide here only
two definitions: for the space $H_{N}^{1,h}\left( \Omega ^{h}\right) $ and
its subspace $H_{N,0}^{1,h}\left( \Omega ^{h}\right) $. Others are similar.
We introduce the space $H_{N}^{1,h}\left( \Omega ^{h}\right) $ as: 
\begin{equation}
\begin{split}
& \hspace{4cm}H_{N}^{1,h}\left( \Omega ^{h}\right) = \\
& =\left\{ 
\begin{array}{c}
Y^{h}(\mathbf{x}^{h}):\left\Vert Y^{h}(\mathbf{x}^{h})\right\Vert
_{H_{N}^{1,h}\left( \Omega ^{h}\right)
}^{2}=\sum\limits_{i,j=1}^{m-1}\int\limits_{a}^{b}\left( Y^{h}\left(
x_{i},y_{j},z\right) \right) ^{2}dz+ \\ 
+\sum\limits_{i,j=1}^{m-1}\int\limits_{a}^{b}\left( Y_{x}^{h}\left(
x_{i},y_{j},z\right) \right)
^{2}dz+\sum\limits_{i,j=1}^{m-1}\int\limits_{a}^{b}\left( Y_{y}^{h}\left(
x_{i},y_{j},z\right) \right) ^{2}dz+ \\ 
+\sum\limits_{i,j=1}^{m-1}\int\limits_{a}^{b}\left( Y_{z}^{h}\left(
x_{i},y_{j},z\right) \right) ^{2}dz<\infty%
\end{array}%
\right\} ,
\end{split}
\label{3.00}
\end{equation}%
\begin{equation}
H_{N,0}^{1,h}\left( \Omega ^{h}\right) =\left\{ Y^{h}(\mathbf{x}^{h})\in
H_{N}^{1,h}\left( \Omega ^{h}\right) :Y^{h}(\mathbf{x}^{h})\mid _{\partial
\Omega ^{h}}=0\right\} .  \label{3.01}
\end{equation}%
By embedding theorem $H_{N}^{1,h}\left( \Omega ^{h}\right) \subset
C_{N}^{h}\left( \overline{\Omega }^{h}\right) $ and 
\begin{equation}
\left\Vert Y^{h}(\mathbf{x}^{h})\right\Vert _{C_{N}^{h}\left( \overline{%
\Omega }^{h}\right) }\leq C\left\Vert Y^{h}(\mathbf{x}^{h})\right\Vert
_{H_{N}^{1,h}\left( \Omega ^{h}\right) },\forall Y^{h}\in H_{N}^{1,h}\left(
\Omega ^{h}\right) ,  \label{3.32}
\end{equation}%
where the number $C=C\left( h_{0},A,\Omega \right) >0$ depends only on
listed parameters, where $h_{0}$ is defined in (\ref{3.26}). Also, it
follows from (\ref{3.26}), (\ref{3.31}) that%
\begin{equation}
\left\Vert Y_{x}^{h}(\mathbf{x}^{h}\right\Vert _{L^{2,h}\left( \Omega
^{h}\right) },\left\Vert Y_{y}^{h}(\mathbf{x}^{h}\right\Vert _{L^{2,h}\left(
\Omega ^{h}\right) }\leq C\left\Vert Y^{h}(\mathbf{x}^{h}\right\Vert
_{L^{2,h}\left( \Omega ^{h}\right) }.  \label{3.320}
\end{equation}%
The following formulas are semidiscrete analogs of (\ref{3.16}): 
\begin{equation}
w^{h}(\mathbf{x}^{h}\mathbf{,}\alpha )=\sum\limits_{n=0}^{N-1}w_{n}^{h}(%
\mathbf{x}^{h})Q_{n}(\alpha ),\text{ }\partial _{\alpha }w^{h}(\mathbf{x}^{h}%
\mathbf{,}\alpha )=\sum\limits_{n=0}^{N-1}w_{n}^{h}(\mathbf{x}%
^{h})Q_{n}^{\prime }(\alpha ).  \label{3.43}
\end{equation}%
Also, let $V^{h}(\mathbf{x}^{h})=\left( w_{0}^{h},\cdots ,w_{N-1}^{h}\right)
^{T}(\mathbf{x}^{h})$. Using (\ref{3.31}) and (\ref{3.43}), we now rewrite
problem (\ref{3.17})-(\ref{3.22}), in the form of partial finite differences
as:%
\begin{equation}
\left. 
\begin{array}{c}
B_{N}V_{z}^{h}\left( \mathbf{x}^{h}\right) +A_{1}^{h}\left( \mathbf{x}%
^{h}\right) V_{x}^{h}\left( \mathbf{x}^{h}\right) +A_{2}^{h}\left( \mathbf{x}%
^{h}\right) V_{y}^{h}\left( \mathbf{x}^{h}\right) + \\ 
+F^{h}\left( V^{h}\left( \mathbf{x}^{h}\right) ,\mathbf{x}^{h}\right)
=0,\quad \mathbf{x}^{h}\in \Omega ^{h},%
\end{array}%
\right.  \label{3.41}
\end{equation}%
\begin{equation}
V^{h}\left( \mathbf{x}^{h}\right) \mid _{\partial \Omega ^{h}}=P^{h}\left( 
\mathbf{x}^{h}\right) .  \label{3.42}
\end{equation}

Suppose that we have found the vector function $V^{h}\left( \mathbf{x}%
^{h}\right) $ satisfying equation (\ref{3.41}) and boundary condition (\ref%
{3.42}). Then it follows from (\ref{3.5}), (\ref{3.8}) and (\ref{3.16}) that
to find the semidiscrete analog $a^{h}\left( \mathbf{x}^{h}\right) $ of the
unknown coefficient $a\left( \mathbf{x}\right) $, we should use: 
\begin{equation}
\left. 
\begin{array}{c}
a^{h}\left( \mathbf{x}^{h}\right) = \\ 
-\left( 1/2d\right) \int\limits_{-d}^{d}\left( \nabla _{\mathbf{x}^{h}}\tau
^{h}/\sqrt{\varepsilon _{r}^{h}}\right) \cdot \nabla _{\mathbf{x}^{h}}\left(
\left( \tau _{z}^{h}(/\sqrt{\varepsilon _{r}^{h}}\right) w^{h}\right) (%
\mathbf{x}^{h},\alpha )d\alpha + \\ 
+\left( 1/2d\right) \int\limits_{-d}^{d}\left( r^{h}(\mathbf{x}^{h},\alpha
)\mu _{s}(\mathbf{x}^{h})\int\limits_{-d}^{d}K(\mathbf{x}^{h},\alpha ,\beta
)\left( r^{h}(\mathbf{x}^{h},\beta )\right) ^{-1}d\beta \right) d\alpha , \\ 
r^{h}(\mathbf{x}^{h},\alpha )=\exp \left( -\left( \sqrt{\varepsilon _{r}}%
/\tau _{z}\right) (\mathbf{x}^{h})\sum\limits_{n=0}^{N-1}w_{n}(\mathbf{x}%
^{h})Q_{n}(\alpha )\right) ,\text{ }\mathbf{x}^{h}\in \Omega ^{h}.%
\end{array}%
\right.  \label{3.44}
\end{equation}

Obviously, the following semidiscrete analog of (\ref{3.1003}) is valid:%
\begin{equation}
\left\{ 
\begin{array}{c}
A_{1}^{h}(\mathbf{x}^{h}),A_{2}^{h}(\mathbf{x}^{h})\in C_{N^{2}}^{h}\left( 
\overline{\Omega }^{h}\right) \text{ and the vector function } \\ 
F^{h}\left( V^{h}\left( \mathbf{x}^{h}\right) ,\mathbf{x}^{h}\right) \text{
is continuously differentiable } \\ 
\text{ with respect to its arguments for }\mathbf{x}^{h}\in \overline{\Omega 
}^{h}.%
\end{array}
\right.  \label{3.440}
\end{equation}%
Let $M^{h}=\max \left( \left\Vert A_{1}^{h}(\mathbf{x}^{h})\right\Vert
_{C_{N^{2}}^{h}},\left\Vert A_{2}^{h}(\mathbf{x}^{h})\right\Vert
_{C_{N^{2}}^{h}}\right) .$ Then%
\begin{equation}
M^{h}\leq M=\max \left( \left\Vert A_{1}(\mathbf{x}^{h})\right\Vert
_{C_{N^{2}}^{h}},\left\Vert A_{2}(\mathbf{x}^{h})\right\Vert
_{C_{N^{2}}^{h}}\right) .  \label{3.441}
\end{equation}%
The following functional $J_{\lambda }^{h}\left( V^{h}\right) $ is the
semidiscrete analog of the functional $J_{\lambda }\left( V\right) $ in (\ref%
{3.240}):

\begin{equation}
\left. J_{\lambda }^{h}\left( V^{h}\right) =\left\Vert \left(
B_{N}V_{z}^{h}+A_{1}^{h}V_{x}^{h}+A_{2}^{h}V_{y}^{h}+F^{h}\left( V^{h}\left( 
\mathbf{x}^{h}\right) ,\mathbf{x}^{h}\right) \right) e^{\lambda
z}\right\Vert _{L_{N}^{2,h}\left( \Omega ^{h}\right) }^{2}.\right.
\label{3.45}
\end{equation}%
Let $R>0$ be an arbitrary number. Define the semidiscrete analog $%
S^{h}\left( R,P^{h}\right) $ of the set $S\left( R,P\right) $ in (\ref{3.23}%
) as:

\begin{equation}
\left. 
\begin{array}{c}
S^{h}\left( R,P^{h}\right) = \\ 
=\left\{ V^{h}\in H_{N}^{1,h}\left( \Omega ^{h}\right) :V^{h}\left( \mathbf{x%
}^{h}\right) \mid _{\partial \Omega ^{h}}=P^{h}\left( \mathbf{x}^{h}\right)
,\left\Vert V^{h}\right\Vert _{H_{N}^{1,h}\left( \Omega ^{h}\right)
}<R\right\} .%
\end{array}%
\right.  \label{3.46}
\end{equation}%
To find an approximate solution $V^{h}\left( \mathbf{x}^{h}\right) $ of
problem (\ref{3.41}), (\ref{3.42}), we consider the following problem:

\textbf{Minimization Problem 2.} \emph{Minimize the functional }$J_{\lambda
}^{h}\left( V^{h}\right) $\emph{\ in (\ref{3.45}) on the set }$\overline{
S^{h}\left( R,P^{h}\right) }$\emph{\ defined in (\ref{3.46}). }

\subsection{Formulations of theorems}

\label{sec:5.2}

\textbf{Lemma 1}. \emph{Consider an\ }$n\times n$\emph{\ matrix }$D$\emph{\
\ \ and assume that the inverse matrix }$D^{-1}$ \emph{exists.\ Then there
exists a number }$\xi =\xi \left( D\right) >0$\emph{\ such that }$\left\Vert
Dx\right\Vert ^{2}\geq \xi \left\Vert x\right\Vert ^{2},\forall x\in \mathbb{%
\ R}^{n},$\emph{\ where }$\left\Vert \cdot \right\Vert $\emph{\ is} \emph{%
the euclidean norm.}

We omit the proof of this lemma since it is well known.

\textbf{Theorem 2.} (Carleman estimate). \emph{Let }$M$\emph{\ be the number
defined in (\ref{3.441}). Assume that (\ref{3.26}) holds.} \emph{\ There
exists a sufficiently large number} $\lambda _{0}=\lambda _{0}(d,M,\Omega
^{h},B_{N}$, $\tau ^{h},\varepsilon _{r}^{h},h_{0})\geq 1$ \emph{depending
only on listed parameters such that the following Carleman estimate holds:}%
\begin{equation}
\left. 
\begin{array}{c}
\left\Vert \left(
B_{N}V_{z}^{h}+A_{1}^{h}V_{x}^{h}+A_{2}^{h}V_{y}^{h}\right) e^{\lambda
z}\right\Vert _{L_{N}^{2,h}\left( \Omega ^{h}\right) }^{2}\geq \\ 
\geq \left( 1/4\right) \cdot \left\Vert \left( B_{N}V_{z}^{h}\right)
e^{\lambda z}\right\Vert _{L_{N}^{2,h}\left( \Omega ^{h}\right) }^{2}+\left(
\lambda ^{2}/8\right) \cdot \left\Vert \left( B_{N}V^{h}\right) e^{\lambda
z}\right\Vert _{L_{N}^{2,h}\left( \Omega ^{h}\right) }^{2}, \\ 
\forall V^{h}\in H_{N,0}^{1,h}\left( \Omega ^{h}\right) ,\forall \lambda
\geq \lambda _{0}.%
\end{array}%
\right.  \label{6.1}
\end{equation}

\textbf{Theorem 3 }(central analytical result). \emph{Assume that (\ref{3.26}%
) holds and let }$S^{h}\left( R,P^{h}\right) $\emph{\ be the set defined in (%
\ref{3.46}). Then: }

1. \emph{At every point }$V^{h}\in \overline{S^{h}\left( R,P^{h}\right) }$%
\emph{\ \ and for every }$\lambda \geq 0$\emph{\ the functional }$J_{\lambda
}^{h}\left( V^{h}\right) $\emph{\ defined in (\ref{3.45}) has the Fr\'{e}
chet derivative }$\left( J_{\lambda }^{h}\right) ^{\prime }\left(
V^{h}\right) \in H_{N,0}^{1,h}\left( \Omega ^{h}\right) .$ \emph{\
Furthermore, the Fr\'{e}chet derivative }$\left( J_{\lambda }^{h}\right)
^{\prime }\left( V^{h}\right) $\emph{\ satisfies the Lipschitz condition
with the number }$\rho >0$\emph{\ is independent on }$V_{1}^{h},V_{2}^{h}:$

\begin{equation}
\begin{split}
& \left\Vert \left( J_{\lambda }^{h}\right) ^{\prime }\left(
V_{2}^{h}\right) -\left( J_{\lambda }^{h}\right) ^{\prime }\left(
V_{1}^{h}\right) \right\Vert _{H_{N}^{1,h}\left( \Omega ^{h}\right) }\leq
\rho \left\Vert V_{2}^{h}-V_{1}^{h}\right\Vert _{H_{N}^{1,h}\left( \Omega
^{h}\right) }, \\
& \hspace{3cm}\forall V_{1}^{h},V_{2}^{h}\in \overline{S^{h}\left(
R,P^{h}\right) }.
\end{split}
\label{6.01}
\end{equation}

2. \emph{There exists a sufficiently large number }$\lambda _{1}$ 
\begin{equation}
\lambda _{1}=\lambda _{1}\left( R,d,M,\Omega ^{h},B_{N},\tau
^{h},\varepsilon _{r}^{h},h_{0}\right) \geq \lambda _{0}\geq 1  \label{6.2}
\end{equation}%
\emph{depending only on listed parameters such that functional (\ref{3.45})
is strictly convex on the set }$\overline{S\left( R,P^{h}\right) },$\emph{\
i.e. there exists a number }$C_{1}=C_{1}\left( R,d,M,\Omega ^{h},B_{N},\tau
^{h},\varepsilon _{r}^{h},h_{0}\right) >0$\emph{\ depending only on listed
parameters such that the following inequality holds:}

\begin{equation}
\left. 
\begin{array}{c}
J_{\lambda }^{h}\left( V_{2}^{h}\right) -J_{\lambda }^{h}\left(
V_{1}^{h}\right) -\left( J_{\lambda }^{h}\right) ^{\prime }\left(
V_{1}^{h}\right) \left( V_{2}^{h}-V_{1}^{h}\right) \geq C_{1}\lambda
^{2}e^{2\lambda a}\left\Vert V_{2}^{h}-V_{1}^{h}\right\Vert
_{H_{N}^{1,h}\left( \Omega ^{h}\right) }^{2}, \\ 
\forall \lambda \geq \lambda _{1},\text{ }\forall V_{1}^{h},V_{2}^{h}\in 
\overline{S^{h}\left( R,P^{h}\right) }.%
\end{array}%
\right.  \label{6.46}
\end{equation}

\emph{3. For each }$\lambda \geq \lambda _{1}$\emph{\ there exists unique
minimizer }$V_{\min ,\lambda }^{h}\in \overline{S^{h}\left( R,P^{h}\right) }$%
\emph{\ of the functional }$J_{\lambda }^{h}\left( V^{h}\right) $\emph{\ on
the set }$\overline{S^{h}\left( R,P^{h}\right) }$\emph{\ and }%
\begin{equation}
\left( J_{\lambda }^{h}\right) ^{\prime }\left( V_{\min ,\lambda
}^{h}\right) \left( V^{h}-V_{\min ,\lambda }^{h}\right) \geq 0,\text{ }
\forall V^{h}\in \overline{S^{h}\left( R,P^{h}\right) }.  \label{6.6}
\end{equation}

\emph{\ }\textbf{Remark 5.1. }\emph{Below }$C_{1}>0$\emph{\ denotes
different numbers depending on the same parameters as ones listed above. }

Let $\delta >0$ be the level of the noise in the data. Our goal now is to
estimate the accuracy of the minimizer $V_{\min ,\lambda }^{h}$ depending on 
$\delta .$ Following the classical concept for ill-posed problems \cite{T},
we assume the existence of the exact solution 
\begin{equation}
V^{h\ast }\in S^{h}\left( R,P^{h\ast }\right)  \label{6.7}
\end{equation}%
of problem (\ref{3.41})-(\ref{3.42}) with the exact, i.e. noiseless data $%
P^{h\ast }.$ Hence, 
\begin{align}
& \hspace{0cm}B_{N}V_{z}^{h\ast }\left( \mathbf{x}^{h}\right)
+A_{1}^{h}\left( \mathbf{x}^{h}\right) V_{x}^{h\ast }\left( \mathbf{x}%
^{h}\right) +A_{2}^{h}\left( \mathbf{x}^{h}\right) V_{x}^{h\ast }\left( 
\mathbf{x}^{h}\right) +  \notag \\
& \hspace{1cm}+F^{h}\left( V^{h\ast }\left( \mathbf{x}^{h}\right) ,\mathbf{x}%
^{h}\right) =0,\mathbf{x}^{h}\in \Omega ^{h},  \label{6.8} \\
& \hspace{2cm}V^{h\ast }\left( \mathbf{x}^{h}\right) \mid _{\partial \Omega
^{h}}=P^{h\ast }\left( \mathbf{x}^{h}\right) .  \label{6.9}
\end{align}%
Let two vector functions $G^{h\ast }\left( \mathbf{x}^{h}\right) $ and $%
G^{h}\left( \mathbf{x}^{h}\right) $ be such that 
\begin{align}
& G^{h\ast }\left( \mathbf{x}^{h}\right) \in S^{h}\left( R,P^{h\ast }\right)
,G^{h}\left( \mathbf{x}^{h}\right) \in S^{h}\left( R,P^{h}\right) ,
\label{6.10} \\
& \hspace{1.5cm}\left\Vert G^{h}-G^{h\ast }\right\Vert _{H_{N}^{1,h}\left(
\Omega ^{h}\right) }<\delta .  \label{6.11}
\end{align}

\textbf{Theorem 4}. \emph{Assume that conditions (\ref{6.7})-(\ref{6.11})
hold. Consider the number }$\lambda _{2},$\emph{\ }%
\begin{equation}
\lambda _{2}=\lambda _{1}\left( 2R,d,M,\Omega ^{h},B_{N},\tau
^{h},\varepsilon _{r}^{h},h_{0}\right) \geq \lambda _{1},  \label{6.12}
\end{equation}%
\emph{where }$\lambda _{1}\left( 2R,d,\Omega ^{h},M^{h},B_{N},\tau
^{h},\varepsilon _{r}^{h},h_{0}\right) $\emph{\ is the number in (\ref{6.2}
). Let }$V_{\min ,\lambda _{2}}^{h}$\emph{\ be the minimizer of functional ( %
\ref{3.45}) on the set }$\overline{S^{h}\left( R,P^{h}\right) },$\emph{\
which was found in Theorem 3. Let }$\alpha \in \left( 0,R\right) $ \emph{\
be a number. Suppose that (\ref{6.7}) is replaced with}%
\begin{equation}
V^{h\ast }\in S^{h}\left( R-\alpha ,P^{h\ast }\right) ,\text{ \emph{where} }%
\alpha >C_{1}\delta .  \label{6.120}
\end{equation}%
\emph{Then the vector function }$V_{\min ,\lambda _{2}}^{h}$\emph{\ belongs
to the open \ set }$S^{h}\left( R,P^{h}\right) $\emph{\ and the following
accuracy estimate holds:}%
\begin{equation*}
\left\Vert V_{\min ,\lambda _{2}}^{h}-V^{h\ast }\right\Vert
_{H_{N}^{1,h}\left( \Omega ^{h}\right) }\leq C_{1}\delta .
\end{equation*}%
\emph{\ }

Consider now the gradient descent method of the minimization of functional (%
\ref{3.45}) on the set $\overline{S^{h}\left( R,P^{h}\right) }.$ Let $%
V_{0}^{h}\in B\left( R/3,P^{h}\right) $ be an arbitrary point of this set.
We take $V_{0}^{h}$ as the starting point of our iterations. Construct the
sequence of the gradient descent method as:%
\begin{equation}
V_{n}^{h}=V_{n-1}^{h}-\beta \left( J_{\lambda _{2}}^{h}\right) ^{\prime
}\left( V_{n-1}^{h}\right) ,n=1,2,...,  \label{6.15}
\end{equation}%
where $\beta >0$ is a small number. Since by Theorem 2 functions $\left(
J_{\lambda _{2}}^{h}\right) ^{\prime }(V_{n-1}^{h}$ $)$ $\in
H_{N,0}^{1,h}\left( \Omega ^{h}\right) $, then all vector functions $%
V_{n}^{h}$ \ have the same boundary conditions $P^{h},$ see (\ref{3.01}) and
(\ref{3.46}).

\textbf{Theorem 5.} \emph{Let conditions of Theorem 4 hold, except that (\ref%
{6.120}) is replaced with} 
\begin{equation*}
V^{h\ast }\in S^{h}\left( \frac{R-\alpha }{3},P^{h\ast }\right) ,\text{
where }\alpha /3>C_{1}\delta .
\end{equation*}%
\emph{Then there exists a sufficiently small number }$\beta >0$\emph{\ and a
number }$\gamma =\gamma \left( \beta \right) \in \left( 0,1\right) $\emph{\
such that in (\ref{6.15}) all functions }$V_{n}^{h}\in S^{h}\left(
R,P^{h}\right) ,$\emph{\ and the following convergence estimates hold}%
\begin{equation}
\left. 
\begin{array}{c}
\left\Vert V_{n}^{h}-V_{\min ,\lambda _{2}}^{h}\right\Vert
_{H_{N}^{1,h}\left( \Omega ^{h}\right) }\leq \beta ^{n}\left\Vert
V_{0}^{h}-V_{\min ,\lambda _{2}}^{h}\right\Vert _{H_{N}^{1,h}\left( \Omega
^{h}\right) }, \\ 
\left\Vert V_{n}^{h}-V^{h\ast }\right\Vert _{H_{N}^{1,h}\left( \Omega
^{h}\right) }\leq C_{1}\delta +\beta ^{n}\left\Vert V_{0}^{h}-V_{\min
,\lambda _{2}}^{h}\right\Vert _{H_{N}^{1,h}\left( \Omega ^{h}\right) }, \\ 
\left\Vert a_{n}^{h}-a^{h\ast }\right\Vert _{L_{N}^{2,h}\left( \Omega
^{h}\right) }\leq C_{1}\delta +\beta ^{n}\left\Vert V_{0}^{h}-V_{\min
,\lambda _{2}}^{h}\right\Vert _{H_{N}^{1,h}\left( \Omega ^{h}\right) },%
\end{array}%
\right.  \label{6.17}
\end{equation}%
\emph{\ where }$a_{n}^{h}\left( \mathbf{x}^{h}\right) $\emph{\ and }$%
a_{n}^{h\ast }\left( \mathbf{x}^{h}\right) $\emph{\ are functions which are
obtained from }$V_{n}^{h}$\emph{\ and }$V^{h\ast }$\emph{\ respectively via
( \ref{3.44}).}

\textbf{Remarks 5.2:}

\begin{enumerate}
\item \emph{By Remark 1.1 estimates (\ref{6.17}) imply that the gradient
descent method (\ref{6.15}) of the minimization of the functional }$%
J_{\lambda }^{h}\left( V^{h}\right) $\emph{\ converges globally for }$%
\lambda =\lambda _{2}$\emph{. Indeed, its starting point }$V_{0}^{h}$ \emph{%
\ \ is an arbitrary point of the set }$S\left( R/3,P^{h}\right) ,$\emph{\ \
and }$R>0$\emph{\ is an arbitrary number.}

\item \emph{We fix }$\lambda =\lambda _{2}$\emph{\ in Theorem 3 only for the
sake of the definiteness. In fact, obvious analogs of these theorems are
valid for any }$\lambda \geq \lambda _{2}.$\emph{\ }

\item \emph{Even though above Theorems 3-5 require sufficiently large values
of the parameter }$\lambda ,$\emph{\ we have numerically established in our
computations in section 7 that actually }$\lambda =5$\emph{\ is sufficient.
A similar observation has been consistently made in all above cited works
about the convexification method. Conceptually, this is similar with the
well known fact from almost any asymptotic theory. Indeed, such a theory
typically claims that if a certain parameter }$X$\emph{\ is sufficiently
large/small, then a certain formula }$Y$\emph{\ is valid with a good
accuracy. However, for any specific numerical implementation with its
specific range of parameters only numerical studies can establish which
exactly value of }$X$\emph{\ is sufficient to obtain a good accuracy of }$Y$ 
\emph{.}

\item \emph{Proofs of Theorems 2, 4 and 5 are similar with the proofs in 
\cite{KTR} of Theorems 4.1, 4.4 and 4.5 respectively. Therefore, we prove in
this paper only Theorem 3.}
\end{enumerate}

\subsection{Proof of Theorem 3}

\label{sec:5.3}

Consider two arbitrary points $V_{1}^{h},V_{2}^{h}\in \overline{S^{h}\left(
R,P^{h}\right) }.$\ Let 
\begin{equation}
W^{h}=V_{2}^{h}-V_{1}^{h}.  \label{7.0}
\end{equation}%
Then by (\ref{3.01}), (\ref{3.46}) and the triangle inequality%
\begin{equation}
W^{h}\in S_{0}^{h}\left( 2R\right) =\left\{ V^{h}\in H_{N,0}^{1,h}\left(
\Omega ^{h}\right) :\left\Vert V^{h}\right\Vert _{H_{N}^{1,h}\left( \Omega
^{h}\right) }\leq 2R\right\} .  \label{7.5}
\end{equation}%
Consider the vector function $F^{h}\left( V^{h}\left( \mathbf{x}^{h}\right) ,%
\mathbf{x}^{h}\right) $ in (\ref{3.41}). It follows from (\ref{3.440}), (\ref%
{7.0}), Remark 5.1 and the multidimensional analog of the Taylor formula 
\cite{V} that the following representation is valid%
\begin{align}
& \hspace{0.1cm}F\left( V_{2}^{h}\left( \mathbf{x}^{h}\right) ,\mathbf{x}%
^{h}\right) =F\left( V_{1}^{h}\left( \mathbf{x}^{h}\right) +W^{h}\left( 
\mathbf{x}^{h}\right) ,\mathbf{x}^{h}\right)  \notag \\
& =F\left( V_{1}^{h}\left( \mathbf{x}^{h}\right) ,\mathbf{x}^{h}\right) +%
\widetilde{F}_{1}\left( V_{1}^{h}\left( \mathbf{x}^{h}\right) ,\mathbf{x}%
^{h}\right) W^{h}\left( \mathbf{x}^{h}\right)  \label{7.6} \\
& \hspace{0.5cm}+\widetilde{F}_{2}\left( V_{1}^{h}\left( \mathbf{x}%
^{h}\right) ,V_{1}^{h}\left( \mathbf{x}^{h}\right) +W^{h}\left( \mathbf{x}%
^{h}\right) ,\mathbf{x}^{h}\right) ,  \notag
\end{align}%
where $\widetilde{F}_{1},\widetilde{F}_{2}$ are such that%
\begin{align}
& \hspace{3cm}\left\vert \widetilde{F}_{1}\left( V_{1}^{h}\left( \mathbf{x}%
^{h}\right) ,\mathbf{x}^{h}\right) \right\vert \leq C_{1},  \label{7.60} \\
& \left\vert \widetilde{F}_{2}\left( V_{1}^{h}\left( \mathbf{x}^{h}\right)
,V_{1}^{h}\left( \mathbf{x}^{h}\right) +W^{h}\left( \mathbf{x}^{h}\right) ,%
\mathbf{x}^{h}\right) \right\vert \leq C_{1}\left( W^{h}\left( \mathbf{x}%
^{h}\right) \right) ^{2}.  \label{7.7}
\end{align}%
In particular, (\ref{7.6}) implies that the expression $\widetilde{F}%
_{1}\left( V_{1}^{h}\left( \mathbf{x}^{h}\right) ,\mathbf{x}^{h}\right)
W^{h}\left( \mathbf{x}^{h}\right) $ is linear with respect to $W^{h}\left( 
\mathbf{x}^{h}\right) .$ By (\ref{7.0}), (\ref{7.6}) and (\ref{7.7})%
\begin{equation}
\left. 
\begin{array}{c}
\left[ L\left( V_{1}^{h}+W^{h}\right) +F\left( V_{1}^{h}\left( \mathbf{x}%
^{h}\right) +W^{h}\left( \mathbf{x}^{h}\right) ,\mathbf{x}^{h}\right) \right]
^{2}= \\ 
=\left[ 
\begin{array}{c}
\left( L\left( V_{1}^{h}\right) +F\left( V_{1}^{h}\left( \mathbf{x}%
^{h}\right) ,\mathbf{x}^{h}\right) \right) + \\ 
+\left( L\left( W^{h}\right) +\widetilde{F}_{1}\left( V_{1}^{h}\left( 
\mathbf{x}^{h}\right) ,\mathbf{x}^{h}\right) W^{h}\left( \mathbf{x}%
^{h}\right) \right) + \\ 
+\widetilde{F}_{2}\left( V_{1}^{h}\left( \mathbf{x}^{h}\right)
,V_{1}^{h}\left( \mathbf{x}^{h}\right) +W^{h}\left( \mathbf{x}^{h}\right) ,%
\mathbf{x}^{h}\right)%
\end{array}%
\right] ^{2}= \\ 
=\left[ L\left( V_{1}^{h}\right) +F\left( V_{1}^{h}\left( \mathbf{x}%
^{h}\right) ,\mathbf{x}^{h}\right) \right] ^{2}+ \\ 
+2\left[ L\left( V_{1}^{h}\right) +F\left( V_{1}^{h}\left( \mathbf{x}%
^{h}\right) ,\mathbf{x}^{h}\right) \right] \left[ L\left( W^{h}\right) +%
\widetilde{F}_{1}\left( V_{1}^{h}\left( \mathbf{x}^{h}\right) ,\mathbf{x}%
^{h}\right) W^{h}\left( \mathbf{x}^{h}\right) \right] + \\ 
+2\left[ L\left( V_{1}^{h}\right) +F\left( V_{1}^{h}\left( \mathbf{x}%
^{h}\right) ,\mathbf{x}^{h}\right) \right] \left[ \widetilde{F}_{2}\left(
V_{1}^{h}\left( \mathbf{x}^{h}\right) ,V_{1}^{h}\left( \mathbf{x}^{h}\right)
+W^{h}\left( \mathbf{x}^{h}\right) ,\mathbf{x}^{h}\right) \right] + \\ 
+\left[ 
\begin{array}{c}
L\left( W^{h}\right) +\widetilde{F}_{1}\left( V_{1}^{h}\left( \mathbf{x}%
^{h}\right) ,\mathbf{x}^{h}\right) W^{h}\left( \mathbf{x}^{h}\right) \\ 
+\widetilde{F}_{2}\left( V_{1}^{h}\left( \mathbf{x}^{h}\right)
,V_{1}^{h}\left( \mathbf{x}^{h}\right) +W^{h}\left( \mathbf{x}^{h}\right) ,%
\mathbf{x}^{h}\right)%
\end{array}%
\right] ^{2}%
\end{array}%
\right.  \label{7.8}
\end{equation}%
Denote%
\begin{equation}
\left. 
\begin{array}{c}
I_{\text{lin}}\left( V_{1}^{h}\mathbf{,}W^{h},\mathbf{x}^{h}\right) = \\ 
=2\left[ L\left( V_{1}^{h}\right) +F\left( V_{1}^{h}\left( \mathbf{x}%
^{h}\right) ,\mathbf{x}^{h}\right) \right] \cdot \left[ L\left( W^{h}\right)
+\widetilde{F}_{1}\left( V_{1}^{h}\left( \mathbf{x}^{h}\right) ,\mathbf{x}%
^{h}\right) W^{h}\left( \mathbf{x}^{h}\right) \right] ,%
\end{array}%
\right.  \label{100}
\end{equation}%
\begin{equation}
\left. 
\begin{array}{c}
I_{\text{nonlin}}^{\left( 1\right) }\left( V_{1}^{h}\mathbf{,}W^{h},\mathbf{x%
}^{h}\right) = \\ 
=2\left[ L\left( V_{1}^{h}\right) +F\left( V_{1}^{h}\left( \mathbf{x}%
^{h}\right) ,\mathbf{x}^{h}\right) \right] \cdot \left[ \widetilde{F}%
_{2}\left( V_{1}^{h}\left( \mathbf{x}^{h}\right) ,V_{1}^{h}\left( \mathbf{x}%
^{h}\right) +W^{h}\left( \mathbf{x}^{h}\right) ,\mathbf{x}^{h}\right) \right]
,%
\end{array}%
\right.  \label{7.10}
\end{equation}%
\begin{equation}
\left. 
\begin{array}{c}
I_{\text{nonlin}}^{\left( 2\right) }\left( V_{1}^{h}\mathbf{,}W^{h},\mathbf{x%
}^{h}\right) = \\ 
\left( L\left( W^{h}\right) +\widetilde{F}_{1}\left( V_{1}^{h}\left( \mathbf{%
x}^{h}\right) ,\mathbf{x}^{h}\right) W^{h}\left( \mathbf{x}^{h}\right) +%
\widetilde{F}_{2}\left( V_{1}^{h}\left( \mathbf{x}^{h}\right)
,V_{1}^{h}\left( \mathbf{x}^{h}\right) +W^{h}\left( \mathbf{x}^{h}\right) ,%
\mathbf{x}^{h}\right) \right) ^{2}%
\end{array}%
\right.  \label{7.11}
\end{equation}%
By (\ref{7.0}) and (\ref{7.8})-(\ref{7.11}) 
\begin{equation}
\begin{split}
& \hspace{0.2cm}\left[ L\left( V_{2}^{h}\right) +F\left( V_{2}^{h}\left( 
\mathbf{x}^{h}\right) ,\mathbf{x}^{h}\right) \right] ^{2}-\left[ L\left(
V_{1}^{h}\right) +F\left( V_{1}^{h}\left( \mathbf{x}^{h}\right) ,\mathbf{x}%
^{h}\right) \right] ^{2}= \\
& =I_{\text{lin}}\left( V_{1}^{h}\mathbf{,}W^{h},\mathbf{x}^{h}\right) +I_{%
\text{nonlin}}^{\left( 1\right) }\left( V_{1}^{h}\mathbf{,}W^{h},\mathbf{x}%
^{h}\right) +I_{\text{nonlin}}^{\left( 2\right) }\left( V_{1}^{h},W^{h},%
\mathbf{x}^{h}\right) .
\end{split}
\label{7.12}
\end{equation}%
It follows from (\ref{7.7}), (\ref{7.10}) and (\ref{7.11}) that%
\begin{equation}
\left. 
\begin{array}{c}
\left\vert I_{\text{nonlin}}^{\left( 1\right) }\left( V_{1}^{h}\mathbf{,}%
W^{h},\mathbf{x}^{h}\right) \right\vert \leq C_{1}\left( W^{h}\left( \mathbf{%
\ \ x}^{h}\right) \right) ^{2},\text{ }\forall W^{h}\in S_{0}^{h}\left(
2R\right) , \\ 
\left\vert I_{\text{nonlin}}^{\left( 2\right) }\left( V_{1}^{h}\mathbf{,}%
W^{h},\mathbf{x}^{h}\right) \right\vert \leq C_{1}\left[ \left(
W_{z}^{h}\left( \mathbf{x}^{h}\right) \right) ^{2}+\left( W^{h}\left( 
\mathbf{x}^{h}\right) \right) ^{2}\right] ,\text{ }\forall W^{h}\in
S_{0}^{h}\left( 2R\right) ,%
\end{array}%
\right.  \label{7.13}
\end{equation}%
where $S_{0}^{h}\left( 2R\right) $ is defined in (\ref{7.5}). By (\ref{3.45}%
), (\ref{7.0}) and (\ref{7.12}) 
\begin{align}
& \hspace{1.8cm}J_{\lambda }^{h}\left( V_{2}^{h}\right) -J_{\lambda
}^{h}\left( V_{1}^{h}\right) =J_{\lambda }^{h}\left( V_{1}^{h}+W^{h}\right)
-J_{\lambda }^{h}\left( V_{1}^{h}\right) =  \label{7.14} \\
& \hspace{0.5cm}=\sum\limits_{n=0}^{N-1}\text{ }\sum\limits_{\left(
i,j\right) =\left( 1,1\right) }^{\left( i,j\right) =\left( m-1,m-1\right)
}\int\limits_{a}^{b}I_{\text{lin }}\left( V_{1}^{h}\left(
x_{i},y_{j},z\right) \mathbf{,}W^{h}\left( x_{i},y_{j},z\right)
,x_{i},y_{j},z\right) e^{2\lambda z}dz  \notag \\
& +\sum\limits_{n=0}^{N-1}\text{ }\sum\limits_{\left( i,j\right) =\left(
1,1\right) }^{\left( i,j\right) =\left( m-1,m-1\right)
}\int\limits_{a}^{b}\sum\limits_{k=1}^{2}I_{\text{nonlin}}^{\left( k\right)
}\left( V_{1}^{h}\left( x_{i},y_{j},z\right) \mathbf{,}W^{h}\left(
x_{i},y_{j},z\right) ,x_{i},y_{j},z\right) e^{2\lambda z}dz.  \notag
\end{align}

Using (\ref{3.32}), (\ref{7.13}) and (\ref{7.14}), we obtain 
\begin{align}
& \left\vert \sum\limits_{n=0}^{N-1}\sum\limits_{\left( i,j\right) =\left(
1,1\right) }^{\left( i,j\right) =\left( m-1,m-1\right)
}\int\limits_{a}^{b}\sum\limits_{k=1}^{2}I_{\text{nonlin}}^{\left( k\right)
}\left( V_{1}^{h}\left( x_{i},y_{j},z\right) \mathbf{,}W^{h}\left(
x_{i},y_{j},z\right) ,x_{i},y_{j},z\right) e^{2\lambda z}dz\right\vert 
\notag \\
& \hspace{2cm}\leq C_{1}e^{2\lambda b}\left\Vert W^{h}\right\Vert
_{H_{N}^{1,h}\left( \Omega ^{h}\right) }^{2},\forall W^{h}\in
S_{0}^{h}\left( 2R\right) .  \label{7.15}
\end{align}%
It follows from (\ref{100})-(\ref{7.11}) that the expression in the second
line of (\ref{7.14}) is linear with respect to $W^{h}$. On the other hand,
the expression in the third line of (\ref{7.14}) is nonlinear with respect
to $W^{h}$.

Consider the linear functional $J_{\lambda ,\text{lin}}^{h}\left(
V_{1}^{h}\right) \left( W^{h}\right) :H_{N,0}^{1,h}\left( \Omega ^{h}\right)
\rightarrow \mathbb{R},$ which is the expression in the second line of (\ref%
{7.14}). It follows from (\ref{3.00})-(\ref{3.320}), (\ref{7.0}), (\ref{7.60}%
), (\ref{100}) and (\ref{7.14}) that 
\begin{equation*}
\left\vert J_{\lambda ,\text{lin}}^{h}\left( V_{1}^{h}\right) \left(
W^{h}\right) \right\vert \leq C_{1}e^{2\lambda b}\left\Vert W^{h}\right\Vert
_{H_{N}^{1,h}\left( \Omega ^{h}\right) },\text{ }\forall W^{h}\in
H_{N}^{1,h}\left( \Omega ^{h}\right) .
\end{equation*}%
Hence, $J_{\lambda ,\text{lin}}^{h}\left( V_{1}^{h}\right) \left(
W^{h}\right) :H_{N,0}^{1,h}\left( \Omega ^{h}\right) \rightarrow \mathbb{R}$
is a bounded linear functional. By Riesz theorem there exists a vector
function $\widetilde{J}_{\lambda ,\text{lin}}^{h}\left( V_{1}^{h}\right) \in
H_{N,0}^{1,h}\left( \Omega ^{h}\right) $ such that 
\begin{equation}
\left( \widetilde{J}_{\lambda ,\text{lin}}^{h}\left( V_{1}^{h}\right)
,Y^{h}\right) =J_{\lambda ,\text{lin}}^{h}\left( V_{1}^{h}\right) \left(
Y^{h}\right) ,\forall Y^{h}\in H_{N,0}^{1,h}\left( \Omega ^{h}\right) .
\label{7.16}
\end{equation}%
Also, using (\ref{7.13})-(\ref{7.16}), we obtain%
\begin{equation}
\lim_{\left\Vert W^{h}\right\Vert _{H_{N}^{1,h}\left( \Omega ^{h}\right)
}\rightarrow 0}\frac{J_{\lambda }^{h}\left( V_{1}^{h}+W^{h}\right)
-J_{\lambda }^{h}\left( V_{1}^{h}\right) -J_{\lambda ,\text{lin}}^{h}\left(
V_{1}^{h}\right) \left( W^{h}\right) }{\left\Vert W^{h}\right\Vert
_{H_{N}^{1,h}\left( \Omega ^{h}\right) }}=0.  \label{7.160}
\end{equation}%
Hence, $J_{\lambda ,\text{lin}}^{h}\left( V_{1}^{h}\right)
:H_{N,0}^{1,h}\left( \Omega ^{h}\right) \rightarrow \mathbb{R}$ is the Fr%
\'{e}chet derivative of the functional $J_{\lambda }^{h}\left( V^{h}\right) $
at the point $V_{1}^{h}.$ We denote it as 
\begin{equation}
\left( J_{\lambda }^{h}\right) ^{\prime }\left( V_{1}^{h}\right)
:=J_{\lambda ,\text{lin}}^{h}\left( V_{1}^{h}\right) .  \label{7.161}
\end{equation}%
The proof of the Lipschitz continuity property (\ref{6.01}) of $\left(
J_{\lambda }^{h}\right) ^{\prime }\left( V^{h}\right) $ is omitted here
since it is completely similar with the proof of Theorem 3.1 of \cite{Bak}.

Using (\ref{7.60}), (\ref{7.7}), (\ref{7.11}) and Cauchy-Schwarz inequality,
we estimate now $I_{\text{nonlin}}^{\left( 2\right) }(V_{1}^{h}$, $W^{h},%
\mathbf{x}^{h})$ from the below,%
\begin{align}
& \hspace{3cm}I_{\text{nonlin}}^{\left( 2\right) }\left( V_{1}^{h}\mathbf{,}%
W^{h},\mathbf{x}^{h}\right) \geq \frac{1}{2}\left( L\left( W^{h}\right)
\right) ^{2}-  \label{7.18} \\
& -\left[ \widetilde{F}_{1}\left( V_{1}^{h}\left( \mathbf{x}^{h}\right) ,%
\mathbf{x}^{h}\right) W^{h}\left( \mathbf{x}^{h}\right) +\widetilde{F}%
_{2}\left( V_{1}^{h}\left( \mathbf{x}^{h}\right) ,V_{1}^{h}\left( \mathbf{x}%
^{h}\right) +W^{h}\left( \mathbf{x}^{h}\right) ,\mathbf{x}^{h}\right) \right]
^{2}\geq  \notag \\
& \hspace{3.1cm}\geq \frac{1}{2}\left( L\left( W^{h}\right) \right)
^{2}-C_{1}\left( W^{h}\left( \mathbf{x}^{h}\right) \right) ^{2}.  \notag
\end{align}%
Thus, Theorem 2, (\ref{7.14}) and (\ref{7.161})-(\ref{7.18}) imply 
\begin{equation}
\begin{split}
& \hspace{0.8cm}J_{\lambda }^{h}\left( V_{1}^{h}+W^{h}\right) -J_{\lambda
}^{h}\left( V_{1}^{h}\right) -\left( J_{\lambda }^{h}\right) ^{\prime
}\left( V_{1}^{h}\right) \left( W^{h}\right) \geq \\
& \hspace{0.5cm}\geq \frac{1}{2}\left\Vert L\left( W^{h}\right) e^{\lambda
z}\right\Vert _{L^{2,h}\left( \Omega ^{h}\right) }^{2}-C_{1}\left\Vert
W^{h}e^{\lambda z}\right\Vert _{L^{2,h}\left( \Omega ^{h}\right) }^{2}\geq \\
& \hspace{0cm}\geq \frac{1}{4}\left\Vert \left( B_{N}W_{z}^{h}\right)
e^{\lambda z}\right\Vert _{L_{N}^{2,h}\left( \Omega ^{h}\right) }^{2}+\frac{%
\lambda ^{2}}{8}\left\Vert \left( B_{N}W^{h}\right) e^{\lambda z}\right\Vert
_{L_{N}^{2,h}\left( \Omega ^{h}\right) }^{2}-C_{1}\left\Vert W^{h}e^{\lambda
z}\right\Vert _{L^{2,h}\left( \Omega ^{h}\right) }^{2}.
\end{split}
\label{7.19}
\end{equation}%
By Lemma 1 there exists a number $\widetilde{C}_{1}=\widetilde{C}_{1}\left(
B_{N},N\right) >0$ such that 
\begin{equation*}
\left\Vert \left( B_{N}W^{h}\right) e^{\lambda z}\right\Vert
_{L_{N}^{2,h}\left( \Omega ^{h}\right) }^{2}\geq \widetilde{C}_{1}\left\Vert
W^{h}e^{\lambda z}\right\Vert _{L_{N}^{2,h}\left( \Omega ^{h}\right) }^{2},%
\text{ }\forall W^{h}\in L^{2,h}\left( \Omega ^{h}\right) ,\forall \lambda
>0,
\end{equation*}%
and the same for $\left\Vert \left( B_{N}W_{z}^{h}\right) e^{\lambda
z}\right\Vert _{L_{N}^{2,h}\left( \Omega ^{h}\right) }^{2}.$ Hence, (\ref%
{7.19}) implies for all $\lambda \geq \lambda _{0}$ 
\begin{align}
& \hspace{2cm}J_{\lambda }^{h}\left( V_{1}^{h}+W^{h}\right) -J_{\lambda
}^{h}\left( V_{1}^{h}\right) -\left( J_{\lambda }^{h}\right) ^{\prime
}\left( V_{1}^{h}\right) \left( W^{h}\right)  \label{7.190} \\
& \geq \widetilde{C}_{1}\left( \left\Vert W_{z}^{h}e^{\lambda z}\right\Vert
_{L_{N}^{2,h}\left( \Omega ^{h}\right) }^{2}+\lambda ^{2}\left\Vert
W^{h}e^{\lambda z}\right\Vert _{L_{N}^{2,h}\left( \Omega ^{h}\right)
}^{2}\right) -C_{1}\left\Vert W^{h}e^{\lambda z}\right\Vert _{L^{2,h}\left(
\Omega ^{h}\right) }^{2},  \notag
\end{align}%
where $\lambda _{0}\geq 1$ was chosen in Theorem 2. Choose the number $%
\lambda _{1}\geq \lambda _{0}$ depending on the parameters listed in (\ref%
{6.2}) such that $\widetilde{C}_{1}\lambda _{1}^{2}/2\geq C_{1}$ and keep in
mind Remark 5.1. Then (\ref{7.190}) implies (\ref{6.46}). Given (\ref{6.46}%
), the existence and uniqueness of the minimizer\ $V_{\min ,\lambda }^{h}\in 
\overline{S^{h}\left( R,P^{h}\right) }$\ of the functional $J_{\lambda
}^{h}\left( V^{h}\right) $\ on the set $\overline{S\left( R,P^{h}\right) }$
for every $\lambda \geq \lambda _{1}$ as well as inequality (\ref{6.6})
follow immediately from a combination of Lemma 2.1 and Theorem 2.1 of \cite%
{Bak}. $\square $

\section{Numerical Studies}

\label{sec:6}

\subsection{Data simulation}

\label{sec:7.1}

We have conducted our numerical studies in the 2d case. Below $\mathbf{x}%
=\left( x,y\right) ,$in (\ref{2.1}) and (\ref{2.2}) $a=1,b=2,A=1/2$ and $%
d=1/2.$ Hence, we obtain 
\begin{equation}
\begin{split}
\Omega & =\left\{ \mathbf{x}:x\in \left( -1/2,1/2\right) ,y\in \left(
1,2\right) \right\} , \\
& \Gamma _{d}=\left\{ \mathbf{x}_{\alpha }=(\alpha ,0):\alpha \in \lbrack
-1/2,1/2]\right\} ,\quad .
\end{split}
\label{5.1}
\end{equation}%
In accordance with the conventional practice in the theory of inverse
problems, we obtain the boundary data (\ref{2.42}) via a computational
simulation, i.e. via the numerical solution of the Forward Problem (\ref%
{2.22}), (\ref{2.23}). Following Theorem 1, we solve this problem via the
solution of the integral equation (\ref{2.32}). To solve this equation, we
consider the partition of the intervals $(1,2)$ and $(-1/2,1/2)$ in (\ref%
{5.1}) as:%
\begin{equation}
\left. 
\begin{array}{c}
1=y_{0}<y_{1}<\cdots <y_{m_{y}}=2,\quad y_{j+1}-y_{j}=h_{y}, \\ 
h_{y}>0,\quad j=0,\cdots ,m_{y}-1, \\ 
-1/2=\alpha _{0}<\alpha _{1}<\cdots <\alpha _{m_{\alpha }}=1/2,\quad \alpha
_{j+1}-\alpha _{j}=h_{\alpha }, \\ 
h_{\alpha }>0,\quad j=0,\cdots ,m_{\alpha }-1,%
\end{array}%
\right.  \label{1}
\end{equation}%
where $m_{y},m_{\alpha }>1$ are two integers. Define the discrete subsets $%
(1,2)_{y}^{h_{y}}$ and $(-1/2,1/2)_{\alpha }^{h_{\alpha }}$ of the intervals 
$(1,2)$ and $(-1/2,1/2)$ as $(1,2)_{y}^{h_{y}}=\left\{ y_{j}\right\}
_{j=0}^{m_{y}}$ and $(-1/2,1/2)_{\alpha }^{h_{\alpha }}=\left\{ \alpha
_{j}\right\} _{j=0}^{m_{\alpha }}$. The fully discrete subset $\Omega
_{discr}^{h}$ of the domain $\Omega $ is:%
\begin{equation}
\begin{split}
& \Omega _{discr}=\left\{ -1/2=x_{0}<x_{1}<\cdots <x_{m}=1/2\right\} \times
(1,2)_{y}^{h_{y}}, \\
& \hspace{1cm}x_{j+1}-x_{j}=h,j=0,\cdots ,m-1,
\end{split}
\label{2}
\end{equation}%
see (\ref{3.25}). Denote the corresponding sets of discrete points:%
\begin{equation}
\left. \mathbf{x}_{discr}=\left\{ (x_{i},y_{k})\in \Omega
_{discr}^{h}\right\} ,\text{ }\mathbf{\alpha }_{discr}=\left\{ (\alpha
_{i},0,\cdots ,0):\alpha _{i}\in (-1/2,1/2)_{\alpha }^{h_{\alpha }}\right\}
.\right.  \label{3.38}
\end{equation}

To compute the numerical solution $u(\mathbf{x}_{discr},\mathbf{\alpha }%
_{discr})$ of the Forward Problem (\ref{2.22}), we need to perform the
numerical integration in the integral equation (\ref{2.32}). We note that
the points in the integrals along the geodesic line $\Gamma (\mathbf{x},%
\mathbf{x}_{\alpha })$ do not necessary belong to the set $\Omega _{discr}$.
Hence, we describe now our numerical interpolation. For any point $%
(x^{\Gamma },y^{\Gamma })\in \Gamma (\mathbf{x},\mathbf{x}_{\alpha })$, we
use the following formula of the numerical interpolation to approximate the
value $U\left( x^{\Gamma },y^{\Gamma }\right) $ of any function $U$ involved
in the numerical computation of the integral over $\Gamma (\mathbf{x},%
\mathbf{x}_{\alpha })$: 
\begin{equation}
\begin{split}
U\left( x^{\Gamma },y^{\Gamma }\right) \approx \frac{1}{hh_{y}}& \left(
x_{j+1}-x^{\Gamma }\right) \left( y_{k+1}-y^{\Gamma }\right) U(x_{j},y_{k})
\\
& +\left( x_{j+1}-x^{\Gamma }\right) \left( y^{\Gamma }-y_{k}\right)
U(x_{j},y_{k+1}) \\
& +\left( x^{\Gamma }-x_{j}\right) \left( y_{k+1}-y^{\Gamma }\right)
U(x_{j+1},y_{k}) \\
& +\left. \left( x^{\Gamma }-x_{j}\right) \left( y^{\Gamma }-y_{k}\right)
U(x_{j+1},y_{k+1})\right] , \\
& \hspace{-2cm}\text{ for }\left( x^{\Gamma },y^{\Gamma }\right) \in \left[
x_{j},x_{j+1}\right] \times \left[ y_{k},y_{k+1}\right] ,\text{ see (\ref{1}
), (\ref{2}).}
\end{split}
\label{3.40}
\end{equation}

As to the kernel $K(\mathbf{x},\alpha ,\beta )$ of the integral operator in (%
\ref{2.22}), we work below with the 2D Henyey-Greenstein function \cite%
{Heino}: 
\begin{equation}
K(\mathbf{x},\alpha ,\beta )=H(\alpha ,\beta )=\frac{1}{2d}\left[ \frac{
1-g^{2}}{1+g^{2}-2g\cos (\alpha -\beta )}\right] ,\quad g=\frac{1}{2}.
\label{5.2}
\end{equation}%
Here, $g=1/2$ means an anisotropic scattering, which is half ballistic with $%
g=0$ an half isotropic scattering with $g=1$ \cite{HT1,HT2,HT3}. We take the
same function $f\left( \mathbf{x}\right) $ as the one in (\ref{2.7}), (\ref%
{2.8}) with $\epsilon =0.05.$

\subsection{Numerical results for the inverse problem}

\label{sec:7.2}

Just as in \cite{KTR}, we set 
\begin{equation}
\mu _{s}(\mathbf{x})=5,\ \mathbf{x}\in \Omega ,\quad \mu _{s}(\mathbf{x}%
)=0,\ \mathbf{x}\in \mathbb{R}^{2}\setminus \Omega .  \label{5.3}
\end{equation}%
We use the third line of (\ref{2.17}) for the coefficient function $a(%
\mathbf{x})$, and we take in this formula 
\begin{equation}
\mu _{a}(\mathbf{x})=\left\{ 
\begin{array}{cc}
c_{a}=const.>0, & \text{inside the tested inclusion,} \\ 
0, & \text{outside the tested inclusion.}%
\end{array}%
\right.  \label{5.4}
\end{equation}%
By (\ref{2.17}), (\ref{5.3}) and (\ref{5.4}) we set: 
\begin{equation}
\text{correct inclusion/background contrast}=1+c_{a}/5,  \label{5.5}
\end{equation}%
\begin{equation}
\begin{split}
& \text{computed inclusion/background contrast} \\
& \hspace{0.5cm}=1+\max \left( \text{computed }\mu _{a}(\mathbf{x})\right)
/5.
\end{split}
\label{5.50}
\end{equation}

In all numerical tests below 
\begin{equation}
\varepsilon _{r}\left( \mathbf{x}\right) =\varepsilon _{r}(x,y)=\left\{ 
\begin{array}{ll}
1+x^{2}\ln (y) & y>1, \\ 
1 & \text{otherwise}.%
\end{array}%
\right.  \label{5.6}
\end{equation}%
This function $\varepsilon _{r}(x,y)$ satisfies conditions (\ref{2.11})-(\ref%
{2.13}). Using the fast marching toolbox "Toolbox Fast Marching" \cite{Peyre}
in MATLAB, we obtain the geodesic lines, and display examples on Figure \ref%
{geodesic_line_c}.

\begin{figure}[htbp]
\centering
\includegraphics[width = 5in]{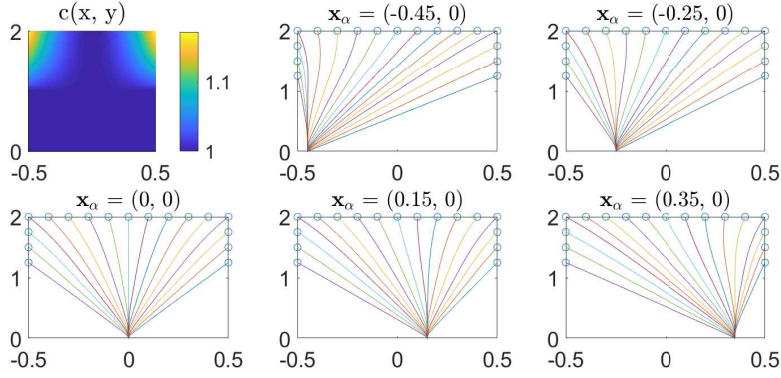}
\caption{Samples of geodesic lines for the function $\protect\varepsilon %
_{r}\left( \mathbf{x}\right) $, which is given in (\protect\ref{5.6}).}
\label{geodesic_line_c}
\end{figure}

The mesh sizes were chosen as $h_{x}=h_{y}=h_{\alpha }=h=1/20$. Hence, we
had total $20\times 20\times N$ unknown parameters in our minimization
procedure. To solve the minimization problem, we have used the Matlab's
built-in function \textbf{fminunc} with the quasi-newton algorithm. The
iterations of the function \textbf{fminunc} were stopped at the iteration
number $k$ as soon as 
\begin{equation*}
\left\vert \nabla J_{\lambda }\left( V_{k}^{h}\right) \right\vert <10^{-2}.
\end{equation*}

The random noise was introduced in the boundary data $g_{1}(\mathbf{x}%
,\alpha )$ in (\ref{3.15}) as: 
\begin{equation}
g_{1}(\mathbf{x},\alpha )=g_{1}(\mathbf{x},\alpha )\left( 1+\delta \cdot
\zeta _{\mathbf{x}}\right) ,\text{ }\mathbf{x}\in \partial \Omega .
\label{5.8}
\end{equation}%
Here $\zeta _{\mathbf{x}}$ is the uniformly distributed random variable in
the interval $[0,1]$ depending on the point $\mathbf{x}\in \partial \Omega $
with $\delta =0.03$ and $\delta =0.05$, which correspond respectively to $%
3\% $ and $5\%$ noise level.

To solve the minimization problem, we need to provide the starting point $%
V_{0}^{h}(\mathbf{x}^{h})$ for iterations. In all numerical tests below we
choose the starting point as the discrete version of the following vector
function $V_{0}(x,y)=\left( w_{0}^{\left( 0\right) },...,w_{N-1}^{\left(
0\right) }\right) ^{T}(x,y):$ 
\begin{equation}
\begin{split}
& w_{n}^{\left( 0\right) }(x,y)=\frac{1}{2}\left( \frac{(A-x)}{2A}
w_{n}(-A,y)+\frac{(x+A)}{2A}w_{n}(A,y)\right) \\
& +\frac{1}{2}\left( \frac{(b-y)}{b-a}w_{n}(x,a)+\frac{(y-a)}{b-a}
w_{n}(x,b)\right) ,\quad n=0,...,N-1.
\end{split}
\label{5.9}
\end{equation}%
Expression (\ref{5.9}) represents the average of linear interpolations of
the boundary condition for $w_{n}\left( x,y\right) $ inside of the square $%
\Omega $ with respect to $x-$direction and $y-$direction.

There are two parameters we need to choose: $N$ and $\lambda $. We find the
optimal pair $\left( N,\lambda \right) =\left( 5,3\right) $ of these
parameters in Test 1, see captions for Figures 2 and 3. Interestingly, the
same optimal pair was found in \cite{KTR} for a similar CIP for the regular
RTE.

\textbf{Remark 6.1}. \emph{To test the computational performance of the
version of the convexification method of this paper, we have chosen
letters-like shapes of abnormalities. This is because letters actually have
complicated shapes for imaging via solutions of CIPs: they are non convex
and have voids. }

We work with the noiseless data in Tests 1-3 and we work with the noisy data
in Test 4.

\textbf{Test 1}. We test the letter `$A$' with $c_{a}=5$ in (\ref{5.4}). We
use this test to figure out optimal values of parameters $N$ and $\lambda $.

First, we select an appropriate value of $N$. We use the value of the norms $%
\left\Vert w_{n}\left( \mathbf{x}\right) \right\Vert _{L_{2}\left( \Omega
\right) }$ to indicate the information contained in $w_{n}\left( \mathbf{x}%
\right) $. Corresponding to the forward problem \eqref{2.22} and \eqref{2.23}
for the case when the functions $\mu _{s}\left( \mathbf{x}\right) $ and $\mu
_{a}\left( \mathbf{x}\right) $ are given in \eqref{5.3} and \eqref{5.4}
respectively, and $c_{a}=5$ in \eqref{5.4}, we calculate norms $\left\Vert
w_{n}\left( \mathbf{x}\right) \right\Vert _{L_{2}\left( \Omega \right) }$
for $n=0,...,11$, and display them in Table \ref{table_w_norm_basis}. One
can see that the $L_{2}\left( \Omega \right) -$norm of the function $w_{n}(%
\mathbf{x})$ decreases very rapidly when the number $n$ is growing, and
these norms, starting from $n=3$ are much less than those for $n=0,1,2$.
More precisely, we have obtained that 
\begin{equation}
\frac{\sum\limits_{n=3}^{11}\left\Vert w_{n}\left( \mathbf{x}\right)
\right\Vert _{L_{2}\left( \Omega \right) }}{\sum\limits_{n=0}^{11}\left\Vert
w_{n}\left( \mathbf{x}\right) \right\Vert _{L_{2}\left( \Omega \right) }}%
=0.0039,  \label{7.1502}
\end{equation}%
which means 0.39\%. We conclude therefore, that we should take in our tests $%
N=3.$

\begin{table}[htbp]
\caption{The $L_{2}\left( \Omega \right) -$norms of functions $w_{n}\left( 
\mathbf{x}\right)$, $n=0,1,...,11$ for the reference Test 1 with $c_{a}=5$
in (\protect\ref{5.4}).}
\label{table_w_norm_basis}\centering
\begin{tabular}{c|c|c|c|c|c|c}
\hline
$n$ & 0 & 1 & 2 & 3 & 4 & 5 \\ \hline
$\left\| w_{n}(\mathbf{x}) \right\|_{L_{2}} $ & 6.5365 & 1.8766 & 0.1924 & 
0.0091 & 0.0071 & 0.0027 \\ \hline
$n$ & 6 & 7 & 8 & 9 & 10 & 11 \\ \hline
$\left\| w_{n}(\mathbf{x}) \right\|_{L_{2}} $ & 0.0057 & 0.0020 & 0.0035 & 
0.0012 & 0.0017 & 0.0008 \\ \hline
\end{tabular}%
\end{table}

\begin{figure}[tbph]
\centering
\includegraphics[width = 4.8in]{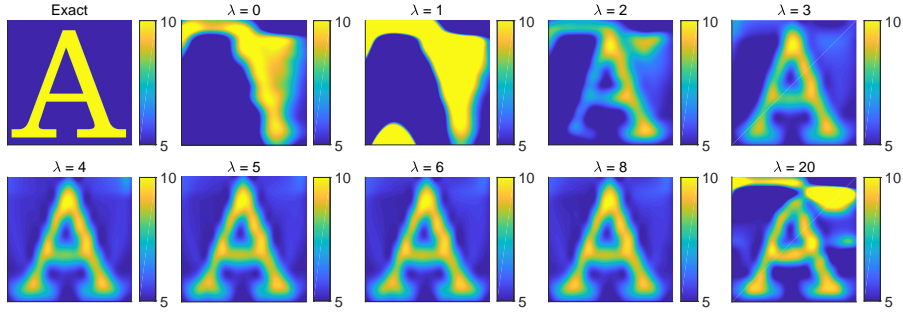}
\caption{Test 1. The reconstructed coefficient $a\left( \mathbf{x}\right) $,
where the function $\protect\mu_{a}\left( \mathbf{x}\right) $ is given in ( 
\protect\ref{5.4}) with $c_{a}=5$ inside of the letter `A'. The goal here is
to test different values of the parameter $\protect\lambda%
=0,1,2,3,4,5,6,8,20 $ for $N=3$. The value of $\protect\lambda $ can be seen
on the top side of each square. The images have a low quality for $\protect%
\lambda =0,1,2,3$. Then the quality is improved and is stabilized at $%
\protect\lambda =5$. Thus, we select $\protect\lambda=5$ as an optimal value
of this parameter for all follow up tests. On the other hand, the last image
is for the case $\protect\lambda =20.$ This image demonstrates that the
quality of the reconstructions deteriorates for too large value of $\protect%
\lambda$.}
\label{plot_NBasis3_diff_Lambda}
\end{figure}

Next, given the value of $N=3$, we select the optimal value of the parameter 
$\lambda $ of the Carleman Weight Function $e^{\lambda z}$ in \eqref{3.45}.
To do this, we test the same letter `A' with $c_{a}=5$ inside of it for
values of the parameter $\lambda =0,1,2,3,4,5,6,8,20.$ Our numerical results
are presented on Figure \ref{plot_NBasis3_diff_Lambda}. We observe that the
images have a low quality for $\lambda =0,1.$ Then the quality is improved,
and it is stabilized at $\lambda =5$. Hence, we treat $\lambda =5$ as the
optimal value of this parameter. Thus, we use in all our tests below%
\begin{equation*}
N=3,\lambda =5.
\end{equation*}

At last, we want to demonstrate numerically again that $N=3$ is indeed a
good choice of $N$ for our optimal value of $\lambda =5.$ Taking $\lambda
=5, $ we test the same letter `A' as above with $c_{a}=5$ in it, but for $%
N=1,2,3,5,7,12.$ The results are displayed in Figure \ref%
{plot_lambda5_diff_NBasis}. One can observe that reconstructions have a low
quality for $N=1,2$. $\ $Next, the reconstructions are basically the same
for $N=3,5,7,12.$ However, the computational cost increases very rapidly
with the increase of $N$. Thus, we conclude that to balance between the
reconstruction accuracy and the computational cost, we should use $N=3$,
which coincides with the above choice.

\begin{figure}[tbph]
\centering
\includegraphics[width = 4.2in]{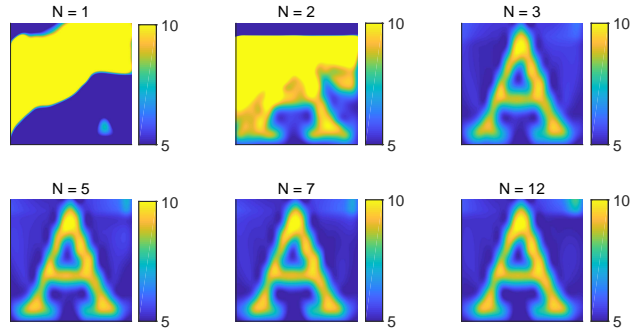}
\caption{Test 1. The reconstructed coefficient $a\left( \mathbf{x}\right) $,
where the function $\protect\mu _{a}\left( \mathbf{x}\right) $ is given in 
\eqref{5.4} with $c_{a}=5$ inside of the letter `A'. We took the optimal
value of the parameter $\protect\lambda =5$ (see Figure \protect\ref%
{plot_NBasis3_diff_Lambda}) and have tested different values of the
parameter $N=1,2,3,5,7,12$. A low quality can be observed for $N=1,2$. The
reconstructions are basically the same for $N=3,5,7,12$. However, the
computational cost increases very rapidly with the increase of $N$. We
conclude, therefore, that to balance between the reconstruction accuracy and
the computational cost, we should use $N=3$. Thus, we use below $\protect%
\lambda =5$ and $N=3$.}
\label{plot_lambda5_diff_NBasis}
\end{figure}

\begin{figure}[tbph]
\centering
\includegraphics[width = 4.5in]{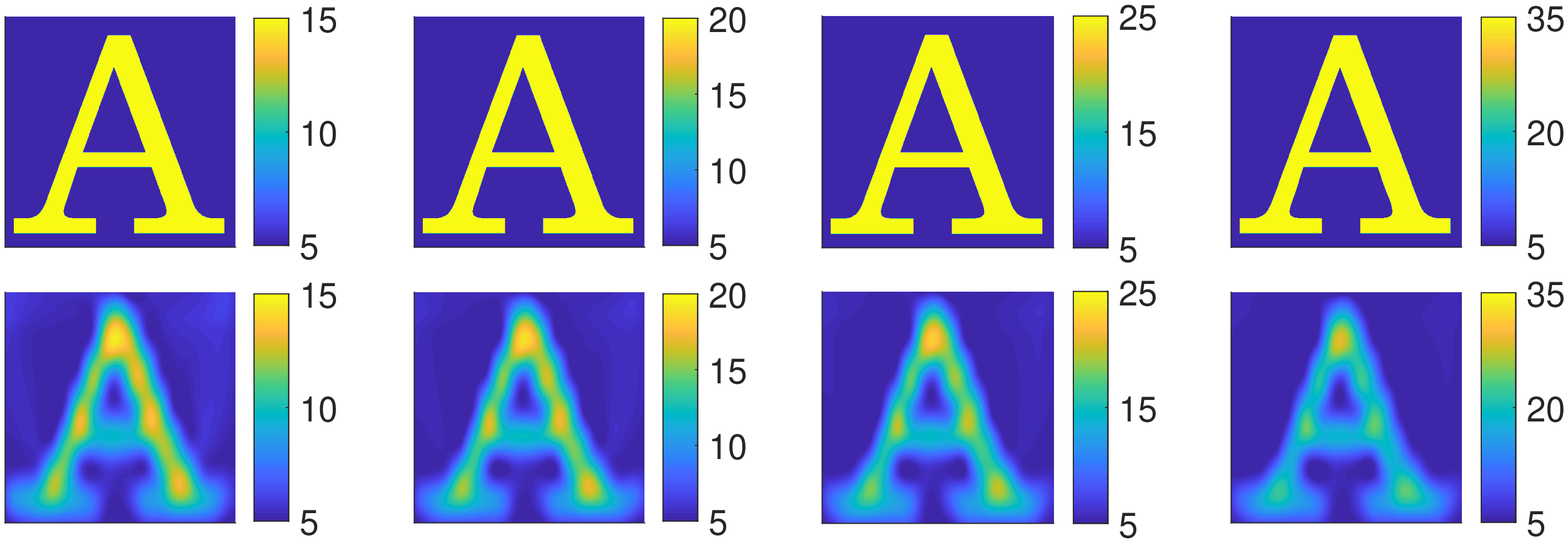}
\caption{Test 2. Exact (top) and reconstructed (bottom) coefficient $a(%
\mathbf{x})$ for $c_{a}=10,15,20,30$ (from left to right ) inside of the
letter `A' as in (\protect\ref{5.4}). Thus, by (\protect\ref{5.5}) the
inclusion/background contrasts now are respectively $3:1$, $4:1$, $5:1$ and $%
6:1$. The image quality remains basically the same for all these values of
the parameter $c_{a}$, although a slight deterioration of this quality can
be observed for $c_{a}=20$ and $c_{a}=30$. The computed inclusion/background
contrasts (\protect\ref{5.5}) are reconstructed accurately.}
\label{plot_A10_A15_A20_A30}
\end{figure}

\textbf{Test 2}. We test the reconstruction of the coefficient $a(\mathbf{x}%
) $ with the shape of the letter `A' where the function $\mu _{a}\left( 
\mathbf{x}\right) $ is given in (\ref{5.4}). We test different values of the
parameter $c_{a}=10,15,20,30$ inside of the letter `A'. Thus, by (\ref{5.5})
the inclusion/background contrasts now are respectively $3:1$, $4:1$, $5:1$
and $6:1$. The function $\varepsilon _{r}(\mathbf{x})=\varepsilon
_{r}^{\left( 1\right) }(\mathbf{x})$ as in (\ref{5.6}). Our computational
results for this test are displayed on Figure \ref{plot_A10_A15_A20_A30}.
One can observe that the quality of these images is good for all four cases,
although it slightly deteriorates for $c_{a}=20$ and $c_{a}=30$. The
computed inclusion/background contrast is accurate, see (\ref{5.50}) and
compare with (\ref{5.5}). 

\begin{figure}[tbph]
\centering
\includegraphics[width = 3.5in]{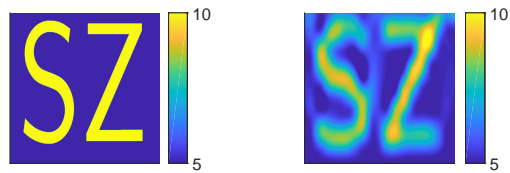}
\caption{Test 3. Exact (left) and reconstructed (right) coefficient $a( 
\mathbf{x})$ for the case when the function $\protect\mu _{a}\left( \mathbf{x%
}\right) $ is given in (\protect\ref{5.4}) with $c_{a}=5$ with the shape of
two letters `SZ'. In (\protect\ref{5.4}) $c_{a}=5$ inside of each of these
two letters and $\protect\mu_{a}\left( \mathbf{x}\right) =0$ outside of each
of these two letters. Here $N=3, \protect\lambda =5$. The quality is good
and the computed inclusion/background contrasts are accurately reconstructed
in both letters, see (\protect\ref{5.5}) and ( \protect\ref{5.50}).}
\label{plot_SZ}
\end{figure}

\textbf{Test 3}. We test the reconstruction of the coefficient $a(\mathbf{x}%
) $ with the shape of two letters `SZ', where the function $\mu _{a}\left( 
\mathbf{x}\right) $ is given in (\ref{5.4}) with $c_{a}=5$ inside of each of
these two letters, and $\mu _{a}\left( \mathbf{x}\right) =0$ outside of each
of these two letters. SZ are two letters in the name of the city (Shenzhen)
were the second and the fifth authors reside. The results are displayed on
Figure \ref{plot_SZ}.

\begin{figure}[tbph]
\centering
\includegraphics[width = 4.5in]{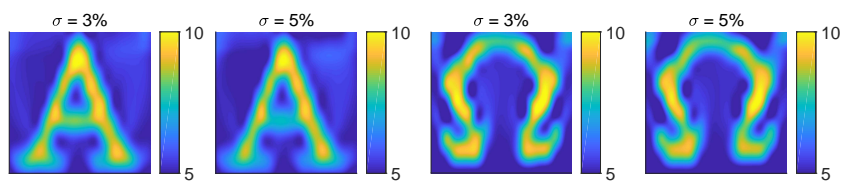}
\caption{Test 4. Reconstructed coefficient $a(\mathbf{x})$ with the shape of
letters `A' and `$\Omega $' with $c_{a}=5$ from noise polluted observation
data as in (\protect\ref{5.8}) with $\protect\delta =0.03$ and $\protect%
\delta =0.05$, i.e. with 3\% and 5\% noise level. One can observe accurate
reconstructions in all four cases. In particular, the inclusion/background
contrasts are reconstructed accurately, see (\protect\ref{5.5}) and (\protect
\ref{5.50}).}
\label{plot_AddNoise}
\end{figure}

\textbf{Test 4}. We now use the noisy data as in (\ref{5.8}) with $\delta
=0.03$ and $\delta =0.05$, i.e. with 3\% and 5\% noise level. We test the
reconstruction of the coefficient $a(\mathbf{x})$ with the shape of either
the letter `A' or the letter `$\Omega $', where the function $\mu _{a}\left( 
\mathbf{x}\right) $ is given in (\ref{5.4}) with $c_{a}=5$ inside of each of
these two letters. The results are displayed on Figure \ref{plot_AddNoise}.
One can observe accurate reconstructions in all four cases. In particular,
the inclusion/background contrasts are reconstructed accurately, see (\ref%
{5.50}) and compare with (\ref{5.5}).

\begin{thebibliography}{10}

\bibitem{B2}
{\sc M.~Asadzadeh and L.~Beilina}, {\em {A stabilized P1 domain decomposition
		finite element method for time harmonic Maxwell's equations}}, Math. Comput.
Simul, 204 (2023), pp.~556--574.

\bibitem{Bak}
{\sc A.~B. Bakushinskii, M.~V. Klibanov, and N.~A. Koshev}, {\em Carleman
	weight functions for a globally convergent numerical method for ill-posed
	{C}auchy problems for some quasilinear {PDE}s}, Nonlinear Anal. Real World
Appl., 34 (2017), pp.~201--224.

\bibitem{BalReview}
{\sc G.~Bal}, {\em Inverse transport theory and applications}, Inverse Probl.,
25 (2009), p.~053001.

\bibitem{Bal1}
{\sc G.~Bal and A.~Jollivet}, {\em Generalized stability estimates in inverse
	transport theory}, Inverse Probl. Imaging, 12 (2018), pp.~59--90.

\bibitem{Bal}
{\sc G.~Bal and A.~Tamasan}, {\em Inverse source problems in transport
	equations}, SIAM J. Math. Anal., 39 (2007), pp.~57--76.

\bibitem{B4}
{\sc L.~Beilina, M.~G. Aram, and E.~M. Karchevskii}, {\em {An adaptive finite
		element method for solving 3D electromagnetic volume integral equation with
		applications in microwave thermometry}}, J. Comput. Phys., 459 (2022),
p.~111122.

\bibitem{B3}
{\sc L.~Beilina and E.~Lindstrom}, {\em {An adaptive finite element/finite
		difference domain decomposition method for applications in microwave
		imaging}}, Electronics, 11 (2022), p.~1359.

\bibitem{BY}
{\sc M.~Bellassoued and M.~Yamamoto}, {\em Carleman estimates and applications
	to inverse problems for hyperbolic systems}, Springer, Japan, 2017.

\bibitem{BW}
{\sc M.~Born and E.~Wolf}, {\em Principles of optics}, Cambridge University
Press, 7th~ed., 1999.

\bibitem{BukhKlib}
{\sc A.~L. Bukhgeim and M.~V. Klibanov}, {\em Uniqueness in the large of a
	class of multidimensional inverse problems}, Soviet Math. Doklady, 17 (1981),
pp.~244--247.

\bibitem{Crist}
{\sc M.~Cristofol, S.~Li, and Y.~Shang}, {\em {Carleman estimates and some
		inverse problems for the coupled quantitative thermoacoustic equations by
		partial boundary layer data. Part II: Some inverse problems}}, Math. Methods
Appl. Sci., published online,  (2023),
\url{https://doi.org/10.1002/mma.9252}.

\bibitem{Fu}
{\sc S.-R. Fu and P.-F. Yao}, {\em Stability in inverse problem of an elastic
	plate with a curved middle surface}, Inverse Probl., 39 (2023), p.~045003.

\bibitem{HT1}
{\sc H.~Fujiwara, K.~Sadiq, and A.~Tamasan}, {\em A {F}ourier approach to the
	inverse source problem in an absorbing and anisotropic scattering medium},
Inverse Probl., 36 (2020), p.~015005.

\bibitem{HT2}
{\sc H.~Fujiwara, K.~Sadiq, and A.~Tamasan}, {\em Numerical reconstruction of
	radiative sources in an absorbing and nondiffusing scattering medium in two
	dimensions}, SIAM J. Imaging Sci., 13 (2020), pp.~535--555.

\bibitem{HT3}
{\sc H.~Fujiwara, K.~Sadiq, and A.~Tamasan}, {\em A source reconstrution method
	in two dimensional radiative transport using boundary data measured on an
	arc}, Inverse Probl., 37 (2021), p.~115005.

\bibitem{Giorgi}
{\sc G.~Giorgi, M.~Brignone, R.~Aramini, and M.~Piana}, {\em {Application of
		the inhomogeneous Lippmann--Schwinger equation to inverse scattering
		problems}}, SIAM J. Appl. Math., 73 (2013), pp.~212--231.

\bibitem{GY}
{\sc F.~G\"{o}lgeleyen and M.~Yamamoto}, {\em Stability for some inverse
	problems for transport equations}, SIAM J. Math. Anal., 48 (2016),
pp.~2319--2344.

\bibitem{Gonch1}
{\sc A.~V. Goncharsky and S.~Y. Romanov}, {\em Iterative methods for solving
	coefficient inverse problems of wave tomography in models with attenuation},
Inverse Probl., 33 (2017), p.~025003.

\bibitem{Gonch2}
{\sc A.~V. Goncharsky and S.~Y. Romanov}, {\em A method of solving the
	coefficient inverse problems of wave tomography}, Comput. Math. Appl., 77
(2019), pp.~967--980.

\bibitem{GN}
{\sc J.~P. Guillement and R.~G. Novikov}, {\em Inversion of weighted {R}adon
	transforms via finite {F}ourier series weight approximation}, Inverse Probl.
Sci. En., 22 (2013), pp.~787--802.

\bibitem{Hassi}
{\sc E.~Hassi, S.-E. Chorfi, and L.~Maniar}, {\em Stable determination of
	coefficients in semilinear parabolic system with dynamic boundary
	conditions}, Inverse Probl., 38 (2022), p.~115007.

\bibitem{Heino}
{\sc J.~Heino, S.~Arridge, J.~Sikora, and E.~Somersalo}, {\em Anisotropic
	effects in highly scattering media}, Phys. Rev. E, 68 (2003), p.~03198.

\bibitem{Is}
{\sc V.~Isakov}, {\em Inverse Problems for Partial Differential Equations},
Springer, New York, 2006.

\bibitem{Kab1}
{\sc S.~I. Kabanikhin, N.~S. Novikov, I.~V. Oseledets, and M.~A. Shishlenin},
{\em Fast toeplitz linear system inversion for solving two-dimensional
	acoustic inverse problem}, J. Inverse Ill-Posed Probl., 23 (2015),
pp.~687--700.

\bibitem{Kab2}
{\sc S.~I. Kabanikhin, K.~K. Sabelfeld, N.~S. Novikov, and M.~A. Shishlenin},
{\em Numerical solution of an inverse problem of coefficient recovering for a
	wave equation by a stochastic projection methods}, Monte Carlo Methods Appl.,
21 (2015), pp.~189--203.

\bibitem{Khoa}
{\sc V.~A. Khoa, M.~V. Klibanov, and L.~H. Nguyen}, {\em Convexification for a
	3{D} inverse scattering problem with the moving point source}, SIAM J. Imag.
Sci., 13 (2020), pp.~871--904.

\bibitem{KY}
{\sc M.~Klibanov and M.~Yamamoto}, {\em Exact controllability for the time
	dependent transport equation}, SIAM J. Control Optim., 46 (2007),
pp.~2071--2195.

\bibitem{Klib92}
{\sc M.~V. Klibanov}, {\em Inverse problems and {C}arleman estimates}, Inverse
Probl., 8 (1992), pp.~575--596.

\bibitem{Klib97}
{\sc M.~V. Klibanov}, {\em Global convexity in a three-dimensional inverse
	acoustic problem}, SIAM J. Math. Anal., 28 (1997), pp.~1371--1388.

\bibitem{Ksurvey}
{\sc M.~V. Klibanov}, {\em Carleman estimates for global uniqueness, stability
	and numerical methods for coefficient inverse problems}, J. Inverse Ill-Posed
Probl., 21 (2013), pp.~477--510.

\bibitem{Klib2017}
{\sc M.~V. Klibanov}, {\em Convexification of restricted {D}irichlet to
	{N}eumann map}, J. Inverse Ill-Posed Probl., 25 (2017), pp.~669--685.

\bibitem{KlibIous}
{\sc M.~V. Klibanov and O.~V. Ioussoupova}, {\em Uniform strict convexity of a
	cost functional for three-dimensional inverse scattering problem}, SIAM J.
Math. Anal, 26 (1995), pp.~147--179.

\bibitem{KL}
{\sc M.~V. Klibanov and J.~Li}, {\em Inverse Problems and Carleman Estimates:
	Global Uniqueness, Global Convergence and Experimental Data}, De Gruyter,
2021.

\bibitem{KTR}
{\sc M.~V. Klibanov, J.~Li, L.~Nguyen, and Z.~Yang}, {\em Convexification
	numerical method for a coefficient inverse problem for the radiative
	transport equation}, SIAM J. Imag. Sci., 16 (2023), pp.~35--63,
\url{https://doi.org/10.1137/22m1509837}.

\bibitem{KTR1}
{\sc M.~V. Klibanov, J.~Li, and Z.~Yang}, {\em Convexification for the
	viscocity solution for a coefficient inverse problem for the radiative
	transfer equation}, arXiv:2302.12474,  (2023).

\bibitem{KP}
{\sc M.~V. Klibanov and S.~Pamyatnykh}, {\em Global uniqueness for a
	coefficient inverse problem for the non-stationary transport equation via
	{C}arleman estimate}, J. Math. Anal. Appl., 343 (2008), pp.~352--365.

\bibitem{KR}
{\sc M.~V. Klibanov and V.~G. Romanov}, {\em A h\"{o}lder stability estimate
	for a coefficient inverse problem for the wave equation with a point source},
Eurasian J. Math. Comp., 10(2) (2022), pp.~11--25.

\bibitem{Lay}
{\sc R.~Y. Lay and Q.~Li}, {\em Parameter reconstruction for general transport
	equation}, SIAM J. Math. Anal., 52 (2020), pp.~2734--2758.

\bibitem{LLM21}
{\sc J.~Li, H.~Liu, and S.~Ma}, {\em Determining a random schr{\" o} dinger
	operator: both potential and source are random}, Commun. Math. Phys., 381
(2021), pp.~527--556.

\bibitem{LLRU15}
{\sc J.~Li, H.~Liu, L.~Rondi, and U.~G.}, {\em {Regularized
		Transformation-Optics Cloaking for the Helmholtz Equation: From Partial Cloak
		to Full Cloak}}, Commun. Math. Phys., 335(2) (2015), pp.~671--712.

\bibitem{McD}
{\sc S.~R. McDowall}, {\em An inverse problem for the transport equation in the
	presence of a {Riemannian} metric}, Pac. J. Math., 216 (2004), pp.~303--326.

\bibitem{Nov}
{\sc R.~Novikov and M.~Santacesaria}, {\em Monochromatic reconstruction
	algorithms for two-dimensional multi-channel inverse problems}, Int. Math.
Res. Not., 2013 (2012), pp.~1205--1229.

\bibitem{Peyre}
{\sc G.~Peyre}, {\em Toolbox fast marching}, MATLAB Central File Exchange,
(2023).

\bibitem{Rom}
{\sc V.~G. Romanov}, {\em Inverse Problems of Mathematical Physics}, VNU Press,
Utrecht, The Netherlands, 1986.

\bibitem{Rom2}
{\sc V.~G. Romanov}, {\em Inverse problems for differential equations with
	memory}, Eurasian J. Math. Comp., 2 (2014), pp.~51--80.

\bibitem{Scales}
{\sc J.~A. Scales, M.~L. Smith, and T.~L. Fischer}, {\em Global optimization
	methods for multimodal inverse problems}, J. Comp. Phys., 103 (1992),
pp.~258--268.

\bibitem{Smirnov}
{\sc A.~V. Smirnov, M.~V. Klibanov, and L.~H. Nguyen}, {\em On an inverse
	source problem for the full radiative transfer equation with incomplete
	data}, SIAM J. Sci. Comput., 41 (2019), pp.~B929--B952.

\bibitem{T}
{\sc A.~N. Tikhonov, A.~V. Goncharsky, V.~V. Stepanov, and A.~G. Yagola}, {\em
	Numerical methods for the solution of ill-posed problems}, Kluwer, London,
1995.

\bibitem{V}
{\sc M.~M. Vajnberg}, {\em Variational method and method of monotone operators
	in the theory of nonlinear equations}, Israel Program for Scientific
Translations, Jerusalem-London, 1973.

\bibitem{Yam}
{\sc M.~Yamamoto}, {\em Carleman estimates for parabolic equations}, Topical
Review. Inverse Probl., 25 (2009), p.~123013.

\end{thebibliography}

\end{document}